\documentclass[12pt]{article}

\usepackage{amsmath}\usepackage{amssymb}\usepackage{amsfonts}\usepackage{amsthm}
\usepackage{color}
\usepackage{graphicx}
\usepackage[english]{babel}\usepackage[latin1]{inputenc}\usepackage{times}\usepackage[T1]{fontenc}
\usepackage{epstopdf}
\usepackage{epsfig}

\newtheorem{thm}{\bf{Theorem}}[section]
\newtheorem{lem}[thm]{\bf{Lemma}}

\newtheorem{rem}[thm]{\bf{Remark}}
\newtheorem{df}[thm]{\bf{Definition}}
\newtheorem{cor}[thm]{\bf{Corollary}}

\newtheorem{fact}[thm]{\bf{Fact}}
\newtheorem{ex}[thm]{\bf{Example}}

\numberwithin{equation}{section}
\usepackage{color}

\newcommand{\dom}{\operatorname{dom}}

\newcommand{\intt}{\operatorname{int}}

\newcommand{\Sin}{\operatorname{Sin}}

\newcommand{\argmin}{\operatornamewithlimits{argmin}}

\title{ Thresholds of Prox-Boundedness of PLQ functions}
\date{\today}
\author{W. Hare\thanks{Mathematics, University of British Columbia, Okanagan Campus (UBCO), 3333 University Way, Kelowna, BC, Canada (UBCO).  Research by this author was supported by NSERC of Canada. {\tt warren.hare@ubc.ca}} \and C. Planiden\thanks{Mathematics, UBCO. Research by this author was supported by UBC UGF and by NSERC of Canada. {\tt chayneplaniden@hotmail.com}}}

\begin{document}

\maketitle\author
\begin{center}\textbf{Dedication}\end{center}\begin{center}
{\em This paper is dedicated to the memory of Jean-Jacques Moreau.}
\vspace{12pt}\end{center}
\setcounter{page}{1}\pagenumbering{arabic}
\begin{abstract}
Introduced in the 1960s, the Moreau envelope has grown to become a key tool in non\-smooth analysis and optimization. Essentially an infimal convolution  with a parametrized norm squared, the Moreau envelope is used in many applications and optimization algorithms. An important aspect in applying the Moreau envelope to nonconvex functions is determining if the function is prox-bounded, that is, if there exists a point $x$ and a parameter $r$ such that the Moreau envelope is finite. The infimum of all such $r$ is called the threshold of prox-boundedness (prox-threshold) of the function $f.$ In this paper, we seek to understand the prox-thresholds of piecewise linear-quadratic (PLQ) functions. (A PLQ function is a function whose domain is a union of finitely many polyhedral sets, and that is linear or quadratic on each piece.) The main result provides a computational technique for determining the prox-threshold for a PLQ function, and further analyzes the behavior of the Moreau envelope of the function using the prox-threshold. We provide several examples to illustrate the techniques and challenges.
\end{abstract}



\section{Introduction}\label{sec:intro}

The \emph{Moreau envelope} $e_r f$ of a proper lower-semicontinuous (lsc) function $f,$ is a smoothing, approximating function that made its first appearance in the mid-1960s \cite{moreau1963,proximite}. It was presented by Jean-Jacques Moreau, together with its associated \emph{proximal mapping} $P_rf,$ as a tool in locating and studying the minima of convex functions. A parametrized function of the \emph{prox-parameter} $r,$ the Moreau envelope is defined as the infimal convolution of $f$ with the scaled norm-squared function $\frac{r}{2}\|\cdot-\bar{x}\|^2.$ It is largely used in matters of minimization of convex functions \cite{howto,smoothing,proximitysums,thresholding,convprox,bregman,fastmoreau,whatshape,rockwets,minmoreau,funcanal}, and more recently it has found a place in non-convex optimization as well \cite{proxlike,proxave,ncproxave,parapr,proxmap,compprox,diffprop,genhess,proxfunc}.

Given a function $f$ and a prox-parameter $r,$ the Moreau envelope is formally defined
$$e_rf(\bar{x}):=\inf\limits_{x\in\dom f}\left\{f(x)+\frac{r}{2}\|x-\bar{x}\|^2\right\}.$$
One of the most inviting properties of the Moreau envelope is that of regularization. Starting with a sufficiently well-behaved function $f,$ such as a convex and lower semicontinuous function, the Moreau envelope is continuously differentiable. In fact, $f$ does not have to be differentiable, or even continuous for that matter, in order for this to happen \cite[Proposition 2.1]{proxmap}. Moreover, the global minimum of $e_rf$ coincides with that of $f,$ in the case where it exists \cite[Proposition 13.37]{rockwets}. So the value of this regularization is clear in matters of minimization of nonsmooth functions.

This paper explores the properties of the \emph{threshold of prox-boundedness} (hereafter referred to simply as \emph{threshold} where convenient). A function $f$ is called \emph{prox-bounded} if there exist $r\geq0$ and $x\in\dom f$ such that $e_rf(x)\in\mathbb{R}.$ The infimum of all such $r$ is called the threshold of prox-boundedness of $f,$ and throughout this paper is denoted by $\bar{r}.$ This number $\bar{r}$ is of interest, as any $r>\bar{r}$ yields $e_rf(x)\in\mathbb{R}$ for all $x$ \cite[Theorem 1.25]{rockwets}, and (if $\bar{r}>0$) any $r$ such that $0\leq r<\bar{r}$ yields $e_rf(x)=-\infty$ for all $x.$ At the threshold itself, the Moreau envelope may be $-\infty$ everywhere, a real number everywhere, or some combination of the two, depending on the characteristics  of $f.$ In this paper we seek to identify the proximal threshold and understand the behavior of the envelope at the threshold.

Thresholds are also of interest due to their importance when dealing with certain programmable tasks in optimization. A prime example is the \emph{proximal point method}, a well-known algorithm used for minimizing functions \cite{martreg,proximite,monops}. The algorithm starts at an arbitrary point $x_0\in\dom f$ and iteratively calculates the proximal mapping
$$x_{i+1}=\argmin\limits_y\left\{f(y)+\frac{r_i}{2}\|y-x_i\|^2\right\}.$$
This method is known to converge to the solution point for convex functions \cite{ontheconv}, and for certain nonconvex functions as well \cite{compprox,ppm,ppa}. There is a question of how to choose the sequence ${r_i},$ and it appears that an ideal starting choice is to use the threshold $\bar{r}$ \cite{dynam}. So for this algorithm, and others that use variants of the proximal point method, it is desirable to be able to calculate the threshold for the function in question. With that in mind, the main result of this work is a computational method of identifying and classifying the thresholds of piecewise linear-quadratic (PLQ) functions.

A PLQ function is a function whose domain is a union of polyhedral sets, and that is linear or quadratic on each of those sets \cite[Definition 10.20]{rockwets} (see Definition \ref{df:plqplq} herein). This is a logical family of functions on which to focus in the present work, as they are commonly used in applications and computational optimization \cite{linesearch,convhullalg,plqmodel,plqopt,boundplq}. They are easily programmable, but complex enough to allow us to illustrate the variety of situations that arise at the threshold.

The organization of this work is as follows. Section \ref{sec:prelim} provides background definitions and presents the method we use to identify the domain of the Moreau envelope, on $\mathbb{R}.$ In Section \ref{sec:quad}, we consider full-domain, quadratic functions on $\mathbb{R}^n.$ In Section \ref{sec:gen} we work with functions that have conic or general polyhedral domains, and we present the main result: computation and classification of the thresholds for PLQ functions. Section \ref{sec:examples} provides several examples that illustrate some special cases and the procedures given in previous sections. Section \ref{sec:conclusion} provides some concluding thoughts, and proposes areas of future research.

\section{Preliminaries}\label{sec:prelim}

\subsection{Notation}
For all that follows, we use the notation $S^n$ for the set of symmetric matrices, $S^n_+$ for the set of symmetric positive-semidefinite matrices, and $S^n_{++}$ for the set of symmetric positive-definite matrices. We introduce the notation $D^n,~D^n_+,\mbox{ and }D^n_{++}$ to represent the sets of diagonal matrices that are arbitrary, positive semidefinite, and positive definite, respectively. For a function $f:\mathbb{R}^n\rightarrow\mathbb{R}\cup\{-\infty,+\infty\},$ we will denote by $\dom f$ the set of points where $f$ is finite, that is,
$$\dom f:=\{x\in\mathbb{R}^n:|f(x)|<+\infty\}.$$
\subsection{Definitions}
\begin{df}\label{df:plqplq}
A function $f:\mathbb{R}^n\rightarrow\mathbb{R}\cup\{+\infty\}$ is called \emph{piecewise linear-quadratic (PLQ)} if $\dom f$ can be represented as the union of finitely many polyhedral sets, relative to each of which $f(x)$ is given by an expression of the form $\frac{1}{2}x^\top Ax+b^\top x+c$ for some scalar $c\in\mathbb{R},$ vector $b\in\mathbb{R}^n,$ and symmetric matrix $A\in S^n.$
\end{df}
\begin{df}\label{df:moreau} The \emph{Moreau envelope} of a proper, lsc function $f:\mathbb{R}^n\rightarrow\mathbb{R}\cup\{+\infty\}$ is denoted $e_rf$ and is defined
$$e_rf(\bar{x}):=\inf\limits_y\left\{f(y)+\frac{r}{2}|y-\bar{x}|^2\right\}.$$
The parameter $r\geq0$ is called the \emph{prox-parameter}, and $x$ is called the \emph{prox-center}.
\end{df}
\begin{df}
A proper lsc function $f:\mathbb{R}^n\rightarrow\mathbb{R}\cup\{+\infty\}$ is \emph{prox-bounded} if there exists $r\geq0$ such that $e_rf(\bar{x})>-\infty$ for some $\bar{x}\in\mathbb{R}^n.$ The infimum of all such $r$ is called the \emph{threshold of prox-boundedness}, and is denoted $\bar{r}.$
\end{df}
For brevity's sake, we refer to the threshold of prox-boundedness of a function simply as its threshold. The goal of this paper is to be able to identify the threshold of any PLQ function, and to describe the behavior of the Moreau envelope at the threshold. We want to be able to say, given any point $\bar{x}\in\mathbb{R}^n,$ whether or not $\bar{x}\in\dom e_{\bar{r}}f.$ It is known that for all $r>\bar{r},$ $\dom e_rf=\mathbb{R}^n,$ and (if $\bar{r}>0$) for any $r\in[0,\bar{r}),$ $\dom e_rf=\emptyset.$ At the threshold itself, however, a variety of situations arise. Depending on the function $f,$ as we see in Examples \ref{ex1}, \ref{ex2}, and \ref{ex3} below, we can have $\dom e_{\bar{r}}f=\mathbb{R}^n,$ $\dom e_{\bar{r}}f=\emptyset,$ or $\emptyset \subsetneq \dom e_{\bar{r}}f \subsetneq \mathbb{R}^n$. We conclude this subsection with a lemma that will be useful in proving some of the results that follow.
\begin{lem}
\label{lem:below}
Let $f:\mathbb{R}^n\rightarrow\mathbb{R}\cup\{+\infty\}$ be proper and lsc. Then $f$ is bounded below if and only if $\bar{r}=0$ and $\dom e_{\bar{r}}f=\mathbb{R}^n.$
\end{lem}
\textbf{Proof:} Notice that
\begin{align*}
f\mbox{ is bounded below}&&&\Leftrightarrow&&\inf_{y\in S}\{f(y)\}>-\infty&\\
&&&\Leftrightarrow&&\inf_{y\in S}\left\{f(y)+\frac{0}{2}\|y-\bar{x}\|^2\right\}>-\infty\mbox{ for all }\bar{x}\in\mathbb{R}^n&\\
&&&\Leftrightarrow&&\bar{r}=0\mbox{ and }\dom e_{\bar{r}}f=\mathbb{R}^n.&\qed
\end{align*}

\subsection{Full-domain single-variable quadratic functions}

We present three examples here, without proof, to show that all three cases above exist in the form of basic functions. The proofs of the example statements are covered by Lemma \ref{lem:quad}. Example \ref{ex2} also demonstrates the importance of the "$\dom e_rf=\mathbb{R}^n$" component of Lemma \ref{lem:below}.
\begin{ex}
\label{ex1}
Let $f(x)=x^2,$ $x\in\mathbb{R}.$ Then $\bar{r}=0$ and $\dom e_{\bar{r}}f=\mathbb{R}.$
\end{ex}
\begin{ex}\label{ex2}
Let $f(x)=x,$ $x\in\mathbb{R}.$ Then $\bar{r}=0$ and $\dom e_{\bar{r}}f=\emptyset.$
\end{ex}
\begin{ex}\label{ex3}
Let $f(x)=-x^2,$ $x\in\mathbb{R}.$ Then $\bar{r}=2$ and $\dom e_{\bar{r}}f=\{0\}.$
\end{ex}
Now we consider a general quadratic function on $\mathbb{R}.$ In the next section, we generalize this result to quadratic functions on $\mathbb{R}^n.$
\begin{lem}
\label{lem:quad}
Let $f:\mathbb{R}\rightarrow\mathbb{R},$ $f(x)=\frac{1}{2}ax^2+bx+c$ be full-domain, i.e., $\dom f=\mathbb{R}.$ Then the threshold of $f$ is $$\bar{r}=\max\{0,-a\},$$ and $\dom e_{\bar{r}}f$ depends on $a$ and $b$ in the following manner.
\begin{itemize}
\item[a)] If $a>0,$ then $\dom e_{\bar{r}}f=\mathbb{R}.$
\item[b)] If $a<0,$ then $\dom e_{\bar{r}}f=\left\{-\frac{b}{a}\right\}.$
\item[c)] If $a=0$ and $b\neq0,$ then $\dom e_{\bar{r}}f=\emptyset.$
\item[d)] If $a=b=0,$ then $\dom e_{\bar{r}}f=\mathbb{R}.$
\end{itemize}
\end{lem}
\textbf{Proof:}
\begin{itemize}
\item[a)] If $a>0,$ then $f$ is bounded below. Hence, $\bar{r}=0$ and $\dom e_{\bar{r}}f=\mathbb{R}$ by Lemma \ref{lem:below}.
\item[b)] If $a<0,$ then for $r\neq-a$ we find the vertex of $\frac{1}{2}ay^2+by+c+\frac{r}{2}(y-x)^2$ by setting the derivative with respect to $y$ equal to $0.$ This gives a critical point $y=\frac{rx-b}{a+r}.$ The second derivative is $a+r,$ so the critical point gives a minimum for all $r>-a,$ and a maximum for all $r<-a.$ Indeed, $r<-a$ results in $\frac{1}{2}ay^2+by+c+\frac{r}{2}(y-x)^2$ being unbounded below. Hence, $\bar{r}=-a.$ Then we evaluate the Moreau envelope at the threshold:
\begin{align*}
e_{-a}f(\bar{x})&=\inf\limits_y\left\{\frac{1}{2}ay^2+by+c+\frac{-a}{2}(y-\bar{x})^2\right\}\\
&=\inf\limits_y\left\{(a\bar{x}+b)y+c-\frac{1}{2}a\bar{x}^2\right\}\\
&=\begin{cases}
c-\frac{b^2}{2a},&\bar{x}=-\frac{b}{a},\\
-\infty,&\mbox{else.}
\end{cases}
\end{align*}
Hence, $\dom e_{\bar{r}}f=\left\{-\frac{b}{a}\right\}.$
\item[c)] If $a=0$ and $b\neq0,$ then for any $r>0$ we have $e_rf(\bar{x})>-\infty$ for all $\bar{x}\in\mathbb{R}.$ This tells us that $\bar{r}=0.$ Then
\begin{align*}
e_{\bar{r}}f(\bar{x})&=\inf\limits_y\{by+c\}\\
&=-\infty\mbox{ for all }\bar{x}\in\mathbb{R}.
\end{align*}
Therefore, $\dom e_{\bar{r}}f=\emptyset.$
\item[d)] If $a=0$ and $b=0,$ then $f$ is constant, and hence bounded below. Lemma \ref{lem:below} applies, and we are done.\qed
\end{itemize}

\section{Full-Domain Quadratic Functions}\label{sec:quad}

Lemma \ref{lem:quad} can be extended to the case $x\in\mathbb{R}^n,$ as we see in Lemma \ref{lem:quad2} and Theorem \ref{thm:quad3}. We begin this section by considering the special case of a quadratic function on $\mathbb{R}^n$ with full domain, whose quadratic coefficient is a diagonal matrix. Recall that we use $D_n,$ $D_n^+,$ and $D_n^{++}$ to denote the sets of $n$-dimensional diagonal, diagonal positive-semidefinite, and diagonal positive definite matrices, respectively.
\begin{lem}
\label{lem:quad2}
Let $f(x)=\frac{1}{2}x^\top Ax+b^\top x+c$ be full-domain, $x\in\mathbb{R}^n,~A\in D^n,~b^\top=(b_1,b_2,\ldots,b_n)\in\mathbb{R}^n,~c\in\mathbb{R}.$ Suppose that (without loss of generality) for $i=1,2,\ldots,n$ the diagonal elements $\lambda_i$ of $A$ are in non-increasing order. Then the threshold of $f$ is $$\bar{r}=\max\{0,-\lambda_n\},$$ and $\dom e_{\bar{r}}f$ depends on $A$ and $b$ in the following manner.
\begin{itemize}
\item[a)] If $A\in D^n_{++},$ then $\dom e_{\bar{r}}f=\mathbb{R}^n.$
\item[b)] If $A\in D^n\setminus D^n_+,$ then $\dom e_{\bar{r}}f=\left\{\bar{x}:\bar{x}_i=-\frac{b_i}{\lambda_i}\mbox{ for all }i\mbox{ such that }\lambda_i=\lambda_n\right\}.$
\item[c)] If $A\in D^n_+\setminus D^n_{++}$ and there exists $i$ such that $\lambda_i=0$ and $b_i\neq0,$ then $\dom e_{\bar{r}}f=\emptyset.$
\item[d)] If $A\in D^n_+\setminus D^n_{++}$ and $b_i=0$ for every $i$ such that $\lambda_i=0,$ then $\dom e_{\bar{r}}f=\mathbb{R}^n.$
\end{itemize}
\end{lem}
\textbf{Proof:} We have
\begin{align}
f(x)&=\frac{1}{2}[x_1,\ldots,x_n]\left[\begin{array}{c c c c}
\lambda_1&0&\cdots&0\\
0&\lambda_2&\cdots&0\\
\vdots&\vdots&\ddots&\vdots\\
0&0&\cdots&\lambda_n
\end{array}\right]\left[\begin{array}{c}
x_1\\
\vdots\\
x_n\end{array}\right]+[b_1,\ldots,b_n]\left[\begin{array}{c}
x_1\\
\vdots\\
x_n\end{array}\right]+c\nonumber\\
&=\frac{1}{2}(\lambda_1x_1^2+\cdots+\lambda_nx_n^2)+(b_1x_1+\cdots+b_nx_n)+c\nonumber\\
&=\left(\frac{\lambda_1}{2}x_1^2+b_1x_1\right)+\left(\frac{\lambda_2}{2}x_2^2+b_2x_2\right)+\cdots+\left(\frac{\lambda_n}{2}x_n^2+b_nx_n\right)+c.
\end{align}
\begin{itemize}
\item[a)] If $A\in D^n_{++},$ then $\lambda_i>0$ for all $i,$ hence, $f$ is bounded below. Therefore, $\bar{r}=0$ and $\dom e_{\bar{r}}f=\mathbb{R}^n$ by Lemma \ref{lem:below}.
\item[b)] If $A\in D^n\setminus D^n_+,$ then $\lambda_n$ is the negative eigenvalue of largest magnitude, since $A$ is ordered. Fix $\bar{x}\in\mathbb{R}^n$ and $r<-\lambda_n,$ and consider the following limit:
\begin{align*}
&\lim\limits_{x_n\to\infty}\left[f(0,\ldots,0,x_n)+\frac{r}{2}|(0,\ldots,0,x_n)-\bar{x}|^2\right]\\
=&\lim\limits_{x_n\to\infty}\left[\frac{\lambda_n}{2}x_n^2+b_nx_n+c+\frac{r}{2}|(0,\ldots,0,x_n)-\bar{x}|^2\right]\\
=&\lim\limits_{x_n\to\infty}\left[\frac{\lambda_n+r}{2}x_n^2+(b_n-r\bar{x}_n)x_n\right]+c+\frac{r}{2}(\bar{x}_1^2+\cdots+\bar{x}_{n-1}^2+\bar{x}_n^2)\\
=&-\infty.
\end{align*}
This gives us that the threshold of $f$ is at least $-\lambda_n.$

Now fix $r>-\lambda_n.$ Then
\begin{align*}
f(x)+\frac{r}{2}|x-\bar{x}|^2&=\frac{1}{2}x^\top Ax+b^\top x+c+\frac{1}{2}(x-\bar{x})^\top(rI)(x-\bar{x})\\
&=\frac{1}{2}x^\top(A+rI)x+b^\top x+\frac{1}{2}\bar{x}^\top(rI)x+c+\bar{x}^\top(rI)\bar{x}.
\end{align*}
Since $r>-\lambda_n,$ then $(A+rI)\in D^n_{++}.$ So $f(x)+\frac{r}{2}|x-\bar{x}|^2$ is strictly convex quadratic, and is therefore bounded below. Hence, $\bar{r}=-\lambda_n.$

Now we consider the Moreau envelope at the threshold:
\begin{align}
e_{\bar{r}}f(\bar{x})&=\inf\limits_y\left\{f(y)-\frac{\lambda_n}{2}|y-\bar{x}|^2\right\}\nonumber\\
&=\inf\limits_y\left\{\frac{1}{2}y^\top Ay+b^\top y+c-\frac{\lambda_n}{2}|y-\bar{x}|^2\right\}\nonumber\\
&=\inf\limits_y\left\{\sum\limits_{i=1}^n\left[\frac{\lambda_i-\lambda_n}{2}y_i^2+(b_i+\lambda_n\bar{x}_i)y_i-\frac{\lambda_n}{2}\bar{x}_i^2\right]+c\right\}.\label{eq:diag}
\end{align}
Notice that $\frac{\lambda_i-\lambda_n}{2}\geq0$ for all $i,$ so that the argument of the infimum above consists of a sum of $n$ single-variable functions, one function of each $y_i,$ that are either strictly convex quadratic (when $\lambda_i>\lambda_n$) or linear (when $\lambda_i=\lambda_n$). In particular, the $n$\textsuperscript{th} such function is linear. Suppose the first $k$ functions are quadratic, and the last $n-k$ functions are linear. Then to find the infimum, we must choose $y_1$ through $y_k$ to be those numbers that give us the vertices of the parabolas $\frac{\lambda_i-\lambda_n}{2}y_i^2+(b_i+\lambda_n\bar{x}_i)y_i-\frac{\lambda_n}{2}\bar{x}_i^2,$ for $i=1,2,\ldots,k.$ That gives us the minimum values for the first $k$ components of the sum in equation (\ref{eq:diag}). For the remaining components, however, we must choose the $y_i$ that give the infima of $(b_i+\lambda_i\bar{x}_i)y_i.$ This means that we will have a finite infimum when $\bar{x}_i=-\frac{b_i}{\lambda_i}$ for each $i=k+1,k+2,\ldots,n,$ but an infimum of $-\infty$ otherwise. Therefore,
\begin{equation}\label{eq:domerf}
\dom e_{\bar{r}}f=\left\{\bar{x}:\bar{x}_i=-\frac{b_i}{\lambda_i},~\lambda_i=\lambda_n\right\}.
\end{equation}
\item[c)] Suppose $A\in D^n_+\setminus D^n_{++},$ and let $k$ be such that $\lambda_k=0$ and $b_k\neq0.$ Fix $\bar{x}\in\mathbb{R}^n$ and consider the Moreau envelope:
$$\inf\limits_y\left\{f(y)+\frac{r}{2}\|y-\bar{x}\|^2\right\}.$$
For any $r>0$ the argument is strictly convex quadratic, so the infimum is a real number. Hence, $\bar{r}=0.$ Now we consider
\begin{align*}
e_{\bar{r}}f(\bar{x})&=\inf\limits_yf(y)\\
&=-\infty\mbox{ for all }\bar{x}\in\mathbb{R}^n,
\end{align*}
since $f$ is linear and non-constant in direction $\bar{x}_k.$ Therefore, $\dom e_{\bar{r}}f=\emptyset.$
\item[d)] Suppose $A\in D^n_+\setminus D^n_{++},$ and $b_i=0$ for all $i$ such that $\lambda_i=0.$ Again we have a finite sum of strictly convex quadratic functions and linear functions, but since $b_i=0$ for every corresponding $\lambda_i=0,$ the linear functions are in fact constant. Hence, the function is bounded below, and we apply Lemma \ref{lem:below} to conclude that $\bar{r}=0$ and $\dom e_{\bar{r}}f=\mathbb{R}^n.$\qed
\end{itemize}

In order to generalize Lemma \ref{lem:quad2} to include all real symmetric matrices, we use the spectral decomposition. Recall that a square matrix $A$ is orthogonally diagonalizable if and only if there exists an orthogonal matrix $Q$ and a diagonal matrix $D$ such that $A=Q^\top DQ.$
\begin{fact}[Fact 8.1.1 \cite{lawa}]
A square matrix $A$ is orthogonally diagonalizable if and only if $A$ is symmetric. Moreover, $D$ is the matrix generated by diagonalizing the vector of eigenvalues of $A.$ This is referred to as the \emph{spectral decomposition} of $A.$
\end{fact}
So if we have a quadratic function $f(x)=\frac{1}{2}x^\top Ax+b^\top x+c$  (where $A$ is symmetric by definition), we are always able to diagonalize $A,$ and the eigenvalues of the resulting diagonal matrix are the same as those of $A.$ The consequence of this is that with a change of variable we will be able to apply Lemma \ref{lem:quad2} to any quadratic, full-domain function. With this tool at our disposal, we present the general form of Lemma \ref{lem:quad2} in Theorem \ref{thm:quad3}.
\begin{thm}\label{thm:quad3}
Let $f(x)=\frac{1}{2}x^\top Ax+b^\top x+c$ be full-domain, $x\in\mathbb{R}^n,$ $A\in S^n,$ $b\in\mathbb{R}^n,$ $c\in\mathbb{R}.$ Let $Q^\top DQ$ be the spectral decomposition of $A,$ and suppose (without loss of generality) that for $i=1,2,\ldots,n$ the diagonal elements $\lambda_i$ of $D$ are in non-increasing order. Then the threshold of $f$ is $$\bar{r}=\max\{0,-\lambda_n\},$$ and $\dom e_{\bar{r}}f$ depends on $D,$ $Q$ and $b$ in the following manner.
\begin{itemize}
\item[a)] If $D\in D^n_{++},$ then $\dom e_{\bar{r}}f=\mathbb{R}^n.$
\item[b)] If $D\in D^n\setminus D^n_+,$ then
\begin{equation}\label{eq:V2}
\dom e_{\bar{r}}f=\left\{\bar{x}:\sum\limits_{j=1}^nq_{ij}\bar{x}_j=-\frac{1}{\lambda_i}\sum\limits_{j=1}^nq_{ij}b_j\mbox{ for all }i\mbox{ with }\lambda_i=\lambda_n\right\}.
\end{equation}
\item[c)] If $D\in D^n_+\setminus D^n_{++}$ and there exists $i$ such that $\lambda_i=0$ and $\sum\limits_{j=1}^nq_{ij}b_j\neq0,$ then $\dom e_{\bar{r}}f=\emptyset.$
\item[d)] If $D\in D^n_+\setminus D^n_{++}$ and $\sum\limits_{j=1}^nq_{ij}b_j=0$ for every $i$ such that $\lambda_i=0,$ then $\dom e_{\bar{r}}f=\mathbb{R}^n.$
\end{itemize}
\end{thm}
\textbf{Proof:} We implement the variable changes $y=Qx$ and $\bar{y}=Q\bar{x}.$ These changes do not affect the threshold, as $Q$ is invertible and, by orthogonality, $Q^{-1}=Q^\top.$ Thus
\begin{align*}
\inf\limits_x\left\{f(x)+\frac{r}{2}|x-\bar{x}|^2\right\}&=\inf\limits_x\left\{f(Q^\top y)+\frac{r}{2}|Q^\top y-Q^\top\bar{y}|^2:y=Qx\right\}\\
&=\inf\limits_y\left\{f(Q^\top y)+\frac{r}{2}|Q^\top y-Q^\top\bar{y}|^2\right\}.
\end{align*}
Further,
\begin{align*}
f(Q^\top y)&=\frac{1}{2}(Q^\top y)^\top A(Q^\top y)+b^\top(Q^\top y)+c\\
&=\frac{1}{2}y^\top QAQ^\top y+(Qb)^\top y+c\\
&=\frac{1}{2}y^\top Dy+(Qb)^\top y+c.
\end{align*}
Now we consider the Moreau envelope,
\begin{align*}
e_rf(Q^\top\bar{y})&=\inf\limits_y\left\{\frac{1}{2}y^\top Dy+(Qb)^\top y+c+\frac{r}{2}\|Q^\top(y-\bar{y})\|^2\right\}\\
&=\inf\limits_y\left\{\frac{1}{2}y^\top Dy+(Qb)^\top y+c+\frac{r}{2}\left[(y-\bar{y})^\top QQ^\top(y-\bar{y})\right]\right\}\\
&=\inf\limits_y\left\{\frac{1}{2}y^\top Dy+(Qb)^\top y+c+\frac{r}{2}\|y-\bar{y}\|^2\right\}.
\end{align*}
Since $D$ is diagonal, we have the form of Lemma \ref{lem:quad2}, with $b$ replaced by $Qb.$ The rest of the proof is analogous to that of Lemma \ref{lem:quad2}.\qed

\begin{rem} An example application of Theorem \ref{thm:quad3} appears in Example \ref{ex:basic}\end{rem}

\section{PLQ Functions}\label{sec:gen}

We next generalize the results we have so far to include functions that have polyhedral domains. We begin by stating some results about the domain of the Moreau envelope; they will be useful in later sections.

\subsection{The Domain of the Moreau Envelope}

In this subsection, we include some useful lemmas about the domain of $e_rf.$ In our first result, we see that the more we restrict the domain of a function, the bigger the domain of the Moreau envelope can be.
\begin{lem}\label{lem:dom1}
Let $f:\dom f\rightarrow\mathbb{R}.$ Suppose $\tilde{f}:\dom\tilde{f}\rightarrow\mathbb{R}$ is such that $\dom\tilde{f}\subseteq\dom f$ and $f(x)=\tilde{f}(x)$ for all $x\in\dom\tilde{f}.$ Then $\dom e_rf\subseteq\dom e_r\tilde{f}.$
\end{lem}
\textbf{Proof:} We have
\begin{align*}
&\inf\limits_{y\in\dom f}\left\{f(y)+\frac{r}{2}\|y-\bar{x}\|^2\right\}>-\infty\mbox{ for all }\bar{x}\in\dom e_rf,\\
\Rightarrow&\inf\limits_{y\in\dom\tilde{f}}\left\{\tilde{f}(y)+\frac{r}{2}\|y-\bar{x}\|^2\right\}>-\infty\mbox{ for all }\bar{x}\in\dom e_rf,
\end{align*}
since $\dom\tilde{f}\subseteq\dom f.$ Therefore, $\dom e_rf\subseteq\dom e_r\tilde{f}.$\qed

Combining Theorem \ref{thm:quad3} with Lemma \ref{lem:dom1}, we have the following corollary.
\begin{cor}\label{cor:dom3}
Let $f:\dom f\rightarrow\subseteq\mathbb{R}^n\rightarrow\mathbb{R},$ $f(x)=\frac{1}{2}x^\top Ax+b^\top x+c$ ($A\in S^n,$ $b\in\mathbb{R}^n,$ $c\in\mathbb{R}$) have threshold $\bar{r}>0.$ For $S\subseteq\dom f,$ let
$$\tilde{f}(x)=\begin{cases}
f(x),&x\in S,\\
\infty,&x\not\in S.
\end{cases}$$ Then
$$\frac{1}{\bar{r}}b\in\dom e_{\bar{r}}f\subseteq\dom e_{\bar{r}}\tilde{f}.$$
\end{cor}
\textbf{Proof:} Using equation (\ref{eq:V2}), we see that substituting $\bar{x}_i=b_i$ satisfies the condition, which gives us that $\frac{1}{\bar{r}}b\in\dom e_{\bar{r}}f.$ Lemma \ref{lem:dom1} completes the proof.\qed

\hbox{}
So for any quadratic function $f$ with $\dom e_{\bar{r}}f\neq\emptyset,$ Corollary \ref{cor:dom3} gives us a point in the domain of the Moreau envelope.

\subsection{Polyhedral Conic Domains}\label{sub:gsc}
Now we are ready to generalize the results of the previous section. We start with a simple case, $f$ quadratic where $\dom f$ is a single, closed, unbounded, conic region. We will change variables to the generalized spherical coordinate form, also known as $n$-spherical coordinates. The variable change is as follows:
$$x\in\mathbb{R}^n\leftrightarrow\begin{array}{r c l}
x_1&=&\rho\cos\phi_1\\
x_2&=&\rho\sin\phi_1\cos\phi_2\\
x_3&=&\rho\sin\phi_1\sin\phi_2\cos\phi_3\\
&&\vdots\\
x_{n-1}&=&\rho\sin\phi_1\cdots\sin\phi_{n-2}\cos\phi_{n-1}\\
x_n&=&\rho\sin\phi_1\cdots\sin\phi_{n-2}\sin\phi_{n-1}
\end{array}
\leftrightarrow 
\begin{array}{l}
(\rho,\phi)\in\mathbb{R}\times\mathbb{R}^{n-1}\\
\rho \geq 0 \\
\phi_1\in[0,2\pi), \\
\phi_2,\phi_3,\ldots,\phi_{n-1}\in[0,\pi].
\end{array}
$$
For ease of notation, we introduce the {\em capital sine-k} function $\Sin_k\phi.$
\begin{df}
Let $\phi=(\phi_1,\phi_2,\ldots,\phi_{n-1}).$ The $\Sin_k$ function is defined
$$\Sin_k\phi:=\prod\limits_{i=1}^k\sin\phi_i.$$
\end{df}
We adopt the conventions $\Sin_0\phi=1$ and $\phi_n=0,$ so that we may write $x_i=\rho\Sin_{i-1}\phi\cos\phi_i$ for all $i=1,2,\ldots,n.$ For a quadratic function $f(x)=\frac{1}{2}x^\top Ax+b^\top x+c,$ the change to $n$-spherical coordinates of the argument of the Moreau envelope results in
\begin{align*}
&\frac{1}{2}x^\top Ax+b^\top x+c+\frac{r}{2}\|x-\bar{x}\|^2\\
=&\frac{1}{2}\sum\limits_{j=1}^n\sum\limits_{i=1}^na_{ij}x_ix_j+\sum\limits_{i=1}^nb_ix_i+c+\frac{r}{2}\sum\limits_{i=1}^n(x_i-\bar{x}_i)^2\\
=&\frac{1}{2}\sum\limits_{j=1}^n\sum\limits_{i=1}^na_{ij}\rho\Sin_{i-1}\phi\cos\phi_i\rho\Sin_{j-1}\phi\cos\phi_j+\sum\limits_{i=1}^nb_i\rho\Sin_{i-1}\phi\cos\phi_i+c\\
&+\frac{r}{2}\sum\limits_{i=1}^n(\rho\Sin_{i-1}\phi\cos\phi_i-\bar{\rho}\Sin_{i-1}\bar{\phi}\cos\bar{\phi}_i)^2\\
=&\rho^2\left(\frac{r}{2}+\frac{1}{2}\sum\limits_{j=1}^n\sum\limits_{i=1}^na_{ij}\Sin_{i-1}\phi\cos\phi_i\Sin_{j-1}\phi\cos\phi_j\right)\\
&+\rho\left(\sum\limits_{i=1}^n(b_i-\bar{\rho}r\Sin_{i-1}\bar{\phi}\cos\bar{\phi}_i)\Sin_{i-1}\phi\cos\phi_i\right)\\
&+c+\frac{\bar{\rho}^2r}{2}\sum\limits_{i=1}^n\Sin_{i-1}^2\bar{\phi}\cos^2\bar{\phi}.
\end{align*}
Define
\begin{align}
G(\phi):=&\sum\limits_{j=1}^n\sum\limits_{i=1}^na_{ij}\Sin_{i-1}\phi\cos\phi_i\Sin_{j-1}\phi\cos\phi_j,\label{eq:gphi}\\ 
H_r(\bar{\rho},\bar{\phi};\phi):=&\sum\limits_{i=1}^n(b_i-\bar{\rho}r\Sin_{i-1}\bar{\phi}\cos\bar{\phi}_i)\Sin_{i-1}\phi\cos\phi_i,\label{eq:h}\\
K_r(\bar{\rho},\bar{\phi}):=&c+\frac{\bar{\rho}^2r}{2}\sum\limits_{i=1}^n\Sin_{i-1}^2\bar{\phi}\cos^2\bar{\phi}_i.\label{eq:k}
\end{align}
Then we have
\begin{equation}\label{eq:ghk}
e_rf(\bar{\rho},\bar{\phi})=\inf\limits_{(\rho,\phi)\in W(S)}\left\{\rho^2\left(\frac{G(\phi)+r}{2}\right)+\rho H_r(\bar{\rho},\bar{\phi};\phi)+K_r(\bar{\rho},\bar{\phi})\right\},
\end{equation}
where $W(x):=(\rho,\phi)$ by the change of variables. Now suppose that $S$ is an unbounded, closed, convex cone. If $S=\mathbb{R}^n,$ then the results of Section \ref{sec:quad} hold. Otherwise, note that $\{\phi:(1,\phi)\in W(S)\}$ is a compact set. Since our expression is quadratic in $\rho,$ it is bounded below if $\frac{1}{2}(G(\phi)+r)>0$ for all $(\rho,\phi)\in W(S),$ and unbounded below if there exists $(\rho,\phi)\in W(S)$ such that $\frac{1}{2}(G(\phi)+r)<0.$ Since $G(\phi)$ is a sum and product of sines and cosines, it is bounded on the compact set $\{\phi:(1,\phi)\in W(S)\},$ and as such it has a minimum. So, defining
\begin{equation}
G:=\min\limits_{(1,\phi)\in S}\{G(\phi)\},\label{eq:g}
\end{equation}
we have
\begin{align*}
&\inf\left\{r:\frac{G(\phi)+r}{2}>0\mbox{ for all }(1,\phi)\in W(S)\right\}\\
=&\inf\{r:r>-G(\phi)\mbox{ for all }(1,\phi)\in W(S)\}\\
=&\inf\{r:r>-G\}\\
=&-G.
\end{align*}
If $G>0,$ then the threshold is $0,$ since it cannot be negative. Hence,
$$\bar{r}=\max\{0,-G\}.$$
Now setting $\bar{r}=\max\{0,-G\},$ we define the following:
\begin{align}
\Phi&:=\{\phi:(1,\phi)\in W(S)\mbox{ and }G(\phi)=G\},\label{eq:Phi}\\
H_{\bar{r}}^+(\bar{\rho},\bar{\phi})&:=\{\phi:\phi\in\Phi\mbox{ and }H_{\bar{r}}(\bar{\rho},\bar{\phi};\phi)\geq0\},\label{eq:hplus}\\
H_{\bar{r}}^{++}(\bar{\rho},\bar{\phi})&:=\{\phi:\phi\in\Phi\mbox{ and }H_{\bar{r}}(\bar{\rho},\bar{\phi};\phi)>0\}.\label{eq:hplusplus}
\end{align}

In the following, recall that a set is said to be polyhedral if it can be expressed as the intersection of a finite number of closed half-spaces \cite[Ex 2.10]{rockwets}.

\begin{thm}\label{thm:polycone}
On $\mathbb{R}^n,$ let $f$ be a quadratic function with $S=\dom f$ a closed, unbounded polyhedral cone. Define $G(\phi),$ $H_r(\bar{\rho},\bar{\phi};\phi),$ $G,$ $\Phi,$ $H_{\bar{r}}^+(\bar{\rho},\bar{\phi}),$ and $H_{\bar{r}}^{++}(\bar{\rho},\bar{\phi})$ as in equations (\ref{eq:gphi}), (\ref{eq:h}), (\ref{eq:g}), (\ref{eq:Phi}), (\ref{eq:hplus}), and (\ref{eq:hplusplus}). Then, using $W(\bar{x})=(\bar{\rho},\bar{\phi}),$ the threshold of $f$ is $$\bar{r}=\max\{0,-G\},$$ and $\dom e_{\bar{r}}f$ depends on $G$ and $H_{\bar{r}}(\bar{\rho},\bar{\phi};\phi)$ in the following manner.
\begin{itemize}
\item[a)] If $G>0,$ then $\dom e_{\bar{r}}f=\mathbb{R}^n.$
\item[b)] If $G \leq 0,$ and $\Phi=H_{\bar{r}}^{++}(\bar{\rho},\bar{\phi}),$ then $\bar{x}\in \dom e_{\bar{r}}f.$
\item[c)] If $G \leq 0,$ and  $\Phi\neq H_{\bar{r}}^+(\bar{\rho},\bar{\phi}),$ then $\bar{x}\notin \dom e_{\bar{r}}f.$
\end{itemize}
\end{thm}
\textbf{Proof:}
\begin{itemize}
\item[a)] If $G>0,$ then $\bar{r}=0$ and we get
\begin{align*}
e_{\bar{r}}f(\bar{\rho},\bar{\phi})&=\inf\limits_{(\rho,\phi)\in W(S)}\{\rho^2G(\phi)+\rho H_{\bar{r}}(\bar{\rho},\bar{\phi};\phi)+K_{\bar{r}}(\bar{\rho},\bar{\phi})\}\\
&>-\infty,
\end{align*}
since $G(\phi)\geq G >0$ for all $(\rho,\phi)\in W(S),$ and hence the argument of the infimum above is a strictly convex (bounded below) function. Therefore, $\dom e_{\bar{r}}f=\mathbb{R}^n.$

\item[b)]
If $G \leq 0,$ then $\bar{r}=-G,$ which gives us that $\rho^2\frac{G(\phi)+\bar{r}}{2}\geq0$ for all $\phi$ with $(1,\phi)\in W(S).$ In fact, $\frac{G(\phi)+\bar{r}}{2}=0$ for all $\phi\in\Phi,$ and $\frac{G(\phi)+\bar{r}}{2}>0$ for all $\phi\not\in\Phi.$ Suppose $(\bar{\rho},\bar{\phi})$ is such that $\Phi=H_{\bar{r}}^{++}(\bar{\rho},\bar{\phi}).$  Consider
    $$\begin{array}{rcl}
    e_{\bar{r}}f(\bar{\rho},\bar{\phi})&=&
    \inf\limits_{(\rho,\phi)\in W(S)} \left\{\rho^2\frac{G(\phi)+\bar{r}}{2} + \rho H_{\bar{r}}(\bar{\rho},\bar{\phi};\phi) + K_{\bar{r}}(\bar{\rho},\bar{\phi})\right\}\\
    &=&
    \min\left\{\begin{array}{l}
    \inf\limits_{\rho \geq 0, \phi\in\Phi}\{\rho^2\frac{G(\phi)+\bar{r}}{2} + \rho H_{\bar{r}}(\bar{\rho},\bar{\phi};\phi) + K_{\bar{r}}(\bar{\rho},\bar{\phi})\}, \\
    \inf\limits_{\rho \geq 0, \phi \notin \Phi} \{\rho^2\frac{G(\phi)+\bar{r}}{2} + \rho H_{\bar{r}}(\bar{\rho},\bar{\phi};\phi) + K_{\bar{r}}(\bar{\rho},\bar{\phi})\}
    \end{array}\right\}.
    \end{array}$$
Working with the first infimum, we note that $\phi\in\Phi$ implies $\frac{G(\phi)+\bar{r}}{2} = 0$, so
    $$\inf\limits_{\rho \geq 0, \phi\in\Phi}\left\{\rho^2\frac{G(\phi)+\bar{r}}{2} + \rho H_{\bar{r}}(\bar{\rho},\bar{\phi};\phi) + K_{\bar{r}}(\bar{\rho},\bar{\phi})\right\}
    = \inf\limits_{\rho \geq 0, \phi\in\Phi}\left\{\rho H_{\bar{r}}(\bar{\rho},\bar{\phi};\phi) + K_{\bar{r}}(\bar{\rho},\bar{\phi})\right\}.$$
As $\Phi=H_{\bar{r}}^{++}(\bar{\rho},\bar{\phi})$, we have $H_{\bar{r}}(\bar{\rho},\bar{\phi};\phi)>0$, so the minimum occurs at $\rho=0$.  That is,
    $$\inf\limits_{\rho \geq 0, \phi\in\Phi}\left\{\rho^2\frac{G(\phi)+\bar{r}}{2} + \rho H_{\bar{r}}(\bar{\rho},\bar{\phi};\phi) + K_{\bar{r}}(\bar{\rho},\bar{\phi})\right\} = K_{\bar{r}}(\bar{\rho},\bar{\phi}).$$
Turning our attention to the second infimum, given any $\phi \notin \Phi$, the inner quadratic is strictly convex. Thus, (using basic calculus) we have\small
    $$\inf_{\rho \geq 0} \left\{\rho^2\frac{G(\phi)+\bar{r}}{2} + \rho H_{\bar{r}}(\bar{\rho},\bar{\phi};\phi) + K_{\bar{r}}(\bar{\rho},\bar{\phi})\right\} =
    \left\{\begin{array}{ll}
    - \frac{(H_{\bar{r}}(\bar{\rho},\bar{\phi};\phi))^2}{2(  G(\phi)+\bar{r})} + K_{\bar{r}}(\bar{\rho},\bar{\phi}) & \mbox{if}~H_{\bar{r}}(\bar{\rho},\bar{\phi};\phi) \leq 0, \\
    K_{\bar{r}}(\bar{\rho},\bar{\phi}) & \mbox{if}~H_{\bar{r}}(\bar{\rho},\bar{\phi};\phi) > 0.
    \end{array}\right.$$\normalsize
Returning to the Moreau envelope calculation, we have that
\begin{equation}\label{eq:erfredone}
    e_{\bar{r}}f(\bar{\rho},\bar{\phi}) =
    \min\left\{\begin{array}{l} K_{\bar{r}}(\bar{\rho},\bar{\phi}),
    \inf\limits_{\stackrel{\phi \notin \Phi}{H_{\bar{r}}(\bar{\rho},\bar{\phi};\phi) \leq 0}} \left\{- \frac{(H_{\bar{r}}(\bar{\rho},\bar{\phi};\phi))^2}{2 ( G(\phi)+\bar{r})} + K_{\bar{r}}(\bar{\rho},\bar{\phi})\right\}
    \end{array}\right\}.
\end{equation}
Since $\Phi=H_{\bar{r}}^{++}(\bar{\rho},\bar{\phi}),$ this simplifies to
    $$e_{\bar{r}}f(\bar{\rho},\bar{\phi}) =
    \min\left\{\begin{array}{l} K_{\bar{r}}(\bar{\rho},\bar{\phi}),
    \inf\limits_{\stackrel{(1,\phi)\in W(S)}{H_{\bar{r}}(\bar{\rho},\bar{\phi};\phi) \leq 0}} \left\{- \frac{(H_{\bar{r}}(\bar{\rho},\bar{\phi};\phi))^2}{2 ( G(\phi)+\bar{r})} + K_{\bar{r}}(\bar{\rho},\bar{\phi})\right\}
    \end{array}\right\}.$$
Finally, noting that $H_{\bar{r}}(\bar{\rho},\bar{\phi};\phi)$ and $G(\phi)$ are continuous functions in $\phi$, and that $\phi$ is bounded, we note that the infimum is over a compact set. Hence, it is obtained:
    $$e_{\bar{r}}f(\bar{\rho},\bar{\phi}) =
    \min\left\{\begin{array}{l} K_{\bar{r}}(\bar{\rho},\bar{\phi}),
    \min\limits_{\stackrel{(1,\phi)\in W(S)}{H_{\bar{r}}(\bar{\rho},\bar{\phi};\phi) \leq 0}} \left\{- \frac{(H_{\bar{r}}(\bar{\rho},\bar{\phi};\phi))^2}{2 ( G(\phi)+\bar{r})} + K_{\bar{r}}(\bar{\rho},\bar{\phi})\right\}
    \end{array}\right\} > - \infty.$$
Therefore, $\bar{x}\in\dom e_{\bar{r}}f.$
\item[c)]
If $(\bar{\rho},\bar{\phi})$ is such that $\Phi \neq H_{\bar{r}}^{+}(\bar{\rho},\bar{\phi})$, then there exists $\hat{\phi}$ such that $\frac{G(\hat{\phi})+\bar{r}}{2}=0$ and $H_{\bar{r}}(\bar{\rho},\bar{\phi};\hat{\phi}) <0$.  Using this, we see that
    $$\begin{array}{rcl}
    e_{\bar{r}}f(\bar{\rho},\bar{\phi})
    &=&
    \inf\limits_{(\rho,\phi)\in W(S)} \left\{\rho^2\frac{G(\phi)+\bar{r}}{2} + \rho H_{\bar{r}}(\bar{\rho},\bar{\phi};\phi) + K_{\bar{r}}(\bar{\rho},\bar{\phi})\right\}\\
    &\leq&
    \inf\limits_{\rho \geq 0}
    \left\{\rho^2\frac{G(\hat{\phi})+\bar{r}}{2} + \rho H_{\bar{r}}(\bar{\rho},\bar{\phi};\hat{\phi}) + K_{\bar{r}}(\bar{\rho},\bar{\phi})\right\}\\
    &=&
    \inf\limits_{\rho \geq 0}
    \{\rho H_{\bar{r}}(\bar{\rho},\bar{\phi};\hat{\phi}) + K_{\bar{r}}(\bar{\rho},\bar{\phi})\} = -\infty.
    \end{array}$$
\qed
\end{itemize}
\begin{rem}
The domain of $e_{\bar{r}}f$ can be identified only in some situations. In particular, the boundary case $G \leq 0 $ and $\Phi=H_{\bar{r}}^+(\bar{\rho},\bar{\phi}) \setminus H_{\bar{r}}^{++}(\bar{\rho},\bar{\phi})$ is not covered by Theorem \ref{thm:polycone}. It is unclear what happens in this situation.
\end{rem}
Before moving to general polyhedral domains, we make one final remark on the domain of the Moreau envelope.

\begin{cor}\label{cor:neq}
On $\mathbb{R}^n,$ let $f$ be a quadratic function with $S=\dom f$ a closed, unbounded polyhedral cone. Define $G(\phi),$ $H_r(\bar{\rho},\bar{\phi};\phi),$ $G,$ $\Phi,$ $H_{\bar{r}}^+(\bar{\rho},\bar{\phi}),$ and $H_{\bar{r}}^{++}(\bar{\rho},\bar{\phi})$ as in equations (\ref{eq:gphi}), (\ref{eq:h}), (\ref{eq:g}), (\ref{eq:Phi}), (\ref{eq:hplus}), and (\ref{eq:hplusplus}). If $G<0,$ then
    $$\dom e_{\bar{r}}f \neq \emptyset \quad \mbox{and} \quad \dom e_{\bar{r}}f \neq \mathbb{R}^n.$$
\end{cor}

\textbf{Proof:} Consider
    $$H_{\bar{r}}(\bar{\rho}, \bar{\phi}; \phi) = \sum\limits_{i=1}^n (b_i-\bar{\rho}\bar{r}\Sin_{i-1}\bar{\phi}\cos\bar{\phi}_i)\Sin_{i-1}\phi\cos\phi_i.$$
We first note that if
    \begin{equation}\label{eq:indom}
    b_i-\bar{\rho}\bar{r}\Sin_{i-1}\bar{\phi}\cos\bar{\phi}_i = 0
    \end{equation}
for all $i$, then $H_{\bar{r}}(\bar{\rho}, \bar{\phi}; \phi) = 0$ for all $\phi$. In this case
    $$\begin{array}{rcl}
    e_{\bar{r}}f(\bar{\rho},\bar{\phi}) &=&
    \inf\limits_{(\rho,\phi)\in W(S)} \{\rho^2\frac{G(\phi)+\bar{r}}{2} + \rho H_{\bar{r}}(\bar{\rho},\bar{\phi};\phi) + K_{\bar{r}}(\bar{\rho},\bar{\phi})\}\\
    &=& \inf\limits_{(\rho,\phi)\in W(S)} \{\rho^2\frac{G(\phi)+\bar{r}}{2} + K_{\bar{r}}(\bar{\rho},\bar{\phi})\}\\
    &\geq& \inf\limits_{(\rho,\phi)\in W(S)} \{K_{\bar{r}}(\bar{\rho},\bar{\phi})\} = K_{\bar{r}}(\bar{\rho},\bar{\phi}),\end{array}$$
as $\rho^2 \frac{G(\phi)+\bar{r}}{2} \geq 0$ for all $(\rho,\phi)\in W(S)$.  Thus, any point $(\bar{\rho},\bar{\phi})$ such that $b_i-\bar{\rho}\bar{r}\Sin_{i-1}\bar{\phi}\cos\bar{\phi}_i = 0$ for all $i$ is in the domain.  Returning equation \eqref{eq:indom} to Cartesian coordinates yields $b_i - \bar{r} \bar{x}_i = 0$, or $\bar{x} = b/\bar{r}$ (so such points clearly exist).

Next, we show that there exists $(\bar{\rho},\bar{\phi})$ such that $H_{\bar{r}}(\bar{\rho}, \bar{\phi}; \phi) < 0$ for some $\phi \in \Phi$.  This means that $(\bar{\rho},\bar{\phi})$ meets the conditions of Theorem \ref{thm:polycone} (c), hence $\dom e_{\bar{r}}f \neq \mathbb{R}^n$.  To see this, select any $\phi \in \Phi$.  Consider the summation
    $$\sum\limits_{i=1}^n (b_i-\bar{\rho}\bar{r}\Sin_{i-1}\bar{\phi}\cos\bar{\phi}_i)\Sin_{i-1}\phi\cos\phi_i.$$

Notice that not all of the factors $\Sin_{i-1}\phi\cos\phi_i$ can be zero. We see this by writing out these terms,
\begin{align*}
\Sin_0\phi\cos\phi_1&=\cos\phi_1,\\
\Sin_1\phi\cos\phi_2&=\sin\phi_1\cos\phi_2,\\
\Sin_2\phi\cos\phi_3&=\sin\phi_1\sin\phi_2\cos\phi_3,\\
&~\vdots\\
\Sin_{n-2}\phi\cos\phi_{n-1}&=\sin\phi_1\cdots\sin\phi_{n-2}\cos\phi_{n-1},\\
\Sin_{n-1}\phi\cos\phi_n&=\sin\phi_1\cdots\sin\phi_{n-1},
\end{align*}
and observing that for the first term to be zero, $\phi_1$ must be either $\frac{\pi}{2}$ or $\frac{3\pi}{2}.$ Then, since $\sin\phi_1=\pm1,$ we must have $\phi_2=\frac{\pi}{2}$ in order for the second term to be zero. Continuing in this manner, we find that $\phi_i=\frac{\pi}{2}$ for all $i\neq1,$ which leaves the last term equal to $\pm1.$ Hence, the summation is never equivalently zero due to $\phi.$

Suppose then, that the $k^{th}$ term, $\Sin_k\phi \cos\phi_k \neq 0$. Selecting $\bar{\phi}_1 = \bar{\phi}_2 = ... \bar{\phi}_{k-1} = \pi/2$ and $\bar{\phi}_k = \bar{\phi}_{k+1} = ... = \bar{\phi}_{n-1} = 0$ yields
    $$\Sin_{i-1}\bar{\phi}\cos\bar{\phi}_i =
    \left\{ \begin{array}{ll}
    0 & \mbox{if}~ i \neq k \\
    1 & \mbox{if}~i=k.
    \end{array}\right.$$
Hence, $b_k-\bar{\rho}\bar{r}\Sin_{k-1}\bar{\phi}\cos\bar{\phi}_k$ can be driven to $-\infty$, while the other terms remain constant, by making $\bar{\rho}$ sufficiently large.  Conversely, selecting $\bar{\phi}_1 = 3\pi/2$, $\bar{\phi}_2 = ... \bar{\phi}_{k-1} = \pi/2$ and $\bar{\phi}_k = \bar{\phi}_{k+1} = ... = \bar{\phi}_{n-1} = 0$ yields
   $$\Sin_{i-1}\bar{\phi}\cos\bar{\phi}_i =
    \left\{ \begin{array}{ll}
    0 & \mbox{if}~ i \neq k \\
    -1 & \mbox{if}~i=k.
    \end{array}\right.$$
Hence, $b_k-\bar{\rho}\bar{r}\Sin_{k-1}\bar{\phi}\cos\bar{\phi}_k$ can be driven to $\infty$, while the other terms remain constant, by making $\bar{\rho}$ sufficiently large.  Therefore, it is always possible to select $(\bar{\rho},\bar{\phi})$ with $H_{\bar{r}}(\bar{\rho}, \bar{\phi}; \phi) < 0$.\qed

\subsection{General polyhedral domains}
Theorem \ref{thm:polycone} covers the case where $\dom f$ is an unbounded polyhedral cone. We now generalize to include all unbounded polyhedral domains. For this, we will need the recession cone, defined as follows.
\begin{df}\cite[Definition 6.33]{rockwets}\label{def:r}
For any point $\bar{x}\in S\subset\mathbb{R}^n,~S\neq\emptyset,$ the \emph{recession cone} $R(\bar{x})$ is the cone defined as
$$R(\bar{x}):=\{x:\bar{x}+\tau x\in S\mbox{ for all }\tau\geq0\}.$$
\end{df}

If $S$ is polyhedral, then $R(\bar{x})$ is the same independent of the choice of $\bar{x}$  \cite[Exercise 6.34]{rockwets}, and we use simply $R.$ If $S$ is bounded, then $R=\{0\}.$ If $S$ is unbounded, then $R$ represents all unbounded directions of $S.$ We will see that in order to understand the threshold, it suffices to focus solely on what happens on $R.$ We first prove that the thresholds themselves are the same on $R$ as on $S,$ in Theorem \ref{thm:rsame} below.
\begin{thm}\label{thm:rsame}
Let $f:S\rightarrow\mathbb{R}$ be a quadratic function with $S$ polyhedral. For any $\hat{x}\in S,$ define $R:=R(\hat{x})+\hat{x}.$ Define
$$\tilde{f}(x):=\begin{cases}
f(x),&x\in R,\\
+\infty,&\mbox{ else.}
\end{cases}$$
Let $\bar{r}_f$ and $\bar{r}_{\tilde{f}}$ be the thresholds of $f$ and $\tilde{f},$ respectively. Then $\bar{r}_{\tilde{f}}=\bar{r}_f.$
\end{thm}
\textbf{Proof:} Let $r>\bar{r}_f.$ Then $\dom e_rf=\mathbb{R}^n,$ so $\dom e_r\tilde{f}=\mathbb{R}^n$ by Lemma \ref{lem:dom1}. This gives us an upper bound on the threshold of $\tilde{f}:$ $\bar{r}_{\tilde{f}}\leq\bar{r}_f.$ Now let $r>\bar{r}_{\tilde{f}}.$ It suffices to show that $\dom e_rf=\mathbb{R}^n,$ since this implies that $r\geq\bar{r}_f.$ Let $\tilde{G}(\phi),$ $\tilde{H}_r(\bar{\rho},\bar{\phi};\phi),$ $\tilde{K}_r(\bar{\rho},\bar{\phi}),$ and $\tilde{G}$ be defined as in equations \eqref{eq:gphi}, \eqref{eq:h}, \eqref{eq:k}, and \eqref{eq:g}, respectively. To see that $\dom e_rf=\mathbb{R}^n,$ suppose that $\bar{x}\not\in\dom e_rf.$ Since $r>\bar{r}_{\tilde{f}},$ we know $\dom e_r\tilde{f}=\mathbb{R}^n,$ so there exists a sequence $\{x_n\}\subseteq S\setminus R$ (where $(\rho_n,\phi_n)=W(x_n)$) such that
\begin{equation}\label{eq:12345}
\lim\limits_{n\rightarrow\infty}\left\{\rho_n^2\frac{\tilde{G}(\phi_n)+r}{2}+\rho_n\tilde{H}_r(\bar{\rho},\bar{\phi};\phi_n)+\tilde{K}_r(\bar{\rho},\bar{\phi})\right\}=-\infty.
\end{equation}
Since $r>\bar{r}_{\tilde{f}},$ we have $\frac{\tilde{G}(\phi)+r}{2}>0$ for all $\phi$ with $(1,\phi)\in W(R).$ Since $\tilde{G}(\phi),$ $\tilde{H}_r(\bar{\rho},\bar{\phi};\phi)$ and $\tilde{K}(\bar{\rho},\bar{\phi})$ are bounded, necessarily $\rho_n\rightarrow\infty.$ By definition of the recession cone, and dropping to a subsequence if necesary, we may assume that $\phi_n\rightarrow\hat{\phi}$ such that $(1,\hat{\phi})\in W(R).$ Since $\frac{\tilde{G}(\hat{\phi})+r}{2}>0$ and $\tilde{G}(\phi)$ is continuous, there exists $N\in\mathbb{N}$ such that $\tilde{G}(\phi_n)+r>\frac{\tilde{G}(\hat{\phi})+r}{2}$ for all $n\geq N.$ This means that\footnotesize
\begin{align*}
\lim\limits_{n\rightarrow\infty}\left(\rho_n^2\frac{\tilde{G}(\phi_n)+r}{2}+\rho_n\tilde{H}_r(\bar{\rho},\bar{\phi};\phi_n)+\tilde{K}_r(\bar{\rho},\bar{\phi})\right)&\geq\lim\limits_{n\rightarrow\infty}\left(\rho_n^2\frac{\tilde{G}(\hat{\phi})+r}{4}+\rho_n\tilde{H}_r(\bar{\rho},\bar{\phi};\phi_n)+\tilde{K}_r(\bar{\rho},\bar{\phi})\right).
\end{align*}\normalsize
Since $\tilde{H}_r(\bar{\rho},\bar{\phi};\phi_n)$ is bounded, say $|\tilde{H}_r(\bar{\rho},\bar{\phi};\phi_n)|\leq L,$ we have that
\begin{align*}
&\lim\limits_{n\rightarrow\infty}\left(\rho_n^2\frac{\tilde{G}(\hat{\phi})+r}{4}+\rho_n\tilde{H}_r(\bar{\rho},\bar{\phi};\phi_n)+\tilde{K}_r(\bar{\rho},\bar{\phi})\right)\\
\geq&\lim\limits_{n\rightarrow\infty}\left(\rho_n^2\frac{\tilde{G}(\hat{\phi})}{4}-\rho_n L+\tilde{K}_r(\bar{\rho},\bar{\phi})\right)\\
=&\infty.
\end{align*}
This is a contradiction to equation (\ref{eq:12345}). Therefore, $\dom e_{\bar{r}}f=\mathbb{R}^n.$\qed

We henceforth drop the subscripts on the threshold and set $\bar{r}_f=\bar{r}_{\tilde{f}}=\bar{r}.$ We now turn our attention to the domain of the Moreau envelope for a polyhedral constrained function.
\begin{thm}\label{thm:recess}
Let $f(x)=\frac{1}{2}x^\top Ax+b^\top x+c,$ ($A\in S^n,$ $b\in\mathbb{R}^n,$ $c\in\mathbb{R}$) be a quadratic function on $S\subseteq\mathbb{R}^n$ with $S$ polyhedral. For any $\hat{x}\in S,$ define $R:=R(\hat{x})+\hat{x}.$ Define
$$\tilde{f}(x):=\begin{cases}
f(x),&x\in R,\\
+\infty,&\mbox{ else.}
\end{cases}$$
Let $\bar{r}$ be the threshold of prox-boundedness of $\tilde{f}.$ For $\tilde{f},$ define $\tilde{G}(\phi),$ $\tilde{H}_r(\bar{\rho},\bar{\phi};\phi),$ $\tilde{G},$ $\tilde{\Phi},$ $\tilde{H}_{\bar{r}}^+(\bar{\rho},\bar{\phi}),$ and $\tilde{H}_{\bar{r}}^{++}(\bar{\rho},\bar{\phi})$ as in equations (\ref{eq:gphi}), (\ref{eq:h}), (\ref{eq:g}), (\ref{eq:Phi}), (\ref{eq:hplus}), and (\ref{eq:hplusplus}). Then the following hold.
\begin{itemize}
\item[a)] If $\tilde{G}>0,$ then $\dom e_{\bar{r}}\tilde{f}=\dom e_{\bar{r}}f=\mathbb{R}^n.$
\item[b)] If $\tilde{G}\leq0,$ and $\phi\in\tilde{\Phi}\Rightarrow(1,\phi)\in\intt R,$ then $\dom e_{\bar{r}}\tilde{f}=\dom e_{\bar{r}}f.$
\item[c)] If $\tilde{G}\leq0$ and $\tilde{\Phi}\neq\tilde{H}_{\bar{r}}^+(\bar{\rho},\bar{\phi}),$ then $\bar{x}\not\in\dom e_{\bar{r}}\tilde{f}$ and $\bar{x}\not\in\dom e_{\bar{r}}f.$
\end{itemize}
\end{thm}
\textbf{Proof:} Notice that the functions $\tilde{G}(\phi)$ and $\tilde{H}_{\bar{r}}(\bar{\rho},\bar{\phi};\phi)$ are the same for $f$ as for $\tilde{f},$ with possibly different domains.
\begin{itemize}
\item[a)] If $\tilde{G}>0,$ then $\bar{r}=0,$ and by the same argument as in the proof of Theorem \ref{thm:polycone} (a) we have $\dom e_{\bar{r}}\tilde{f}=\mathbb{R}^n.$ Suppose $\dom e_{\bar{r}}f\neq\mathbb{R}^n.$ Then there exists $(\bar{\rho},\bar{\phi})$ such that $e_{\bar{r}}f(\bar{\rho},\bar{\phi})=-\infty.$ That is,
\begin{equation}\label{eq:recessionsame}
\inf\limits_{(\rho,\phi)\in W(S\setminus R)}\left\{\rho^2\frac{\tilde{G}(\phi)}{2}+\rho\tilde{H}_{\bar{r}}(\bar{\rho},\bar{\phi};\phi)+\tilde{K}_{\bar{r}}(\bar{\rho},\bar{\phi})\right\}=-\infty.
\end{equation}
In order for equation (\ref{eq:recessionsame}) to be true, we must have a sequence $\{(\rho_n,\phi_n)\}_{n=1}^\infty\subseteq W(S\setminus R)$ such that
\begin{equation}\label{eq:limit}
\lim\limits_{n\rightarrow\infty}\left(\rho_n^2\frac{\tilde{G}(\phi_n)}{2}+\rho_n\tilde{H}_{\bar{r}}(\bar{\rho},\bar{\phi};\phi_n)+\tilde{K}_{\bar{r}}(\bar{\rho},\bar{\phi})\right)=-\infty.
\end{equation}
As $\tilde{G}(\phi)>0$ for all $\phi$ with $(1,\phi)\in W(R),$ by the same argument as in the proof of Theorem \ref{thm:rsame}, we get a contradiction to equation \eqref{eq:recessionsame} and we conclude that $\dom e_{\bar{r}}f=\mathbb{R}^n.$

\item[b)] By Lemma \ref{lem:dom1}, we have $\dom e_{\bar{r}}f\subseteq\dom e_{\bar{r}}\tilde{f}.$ Suppose there exists $\bar{x}\in\dom e_{\bar{r}}\tilde{f}\setminus\dom e_{\bar{r}}f.$ As in part (a), this implies that we have a sequence $\{(\rho_n,\phi_n)\}_{n=1}^\infty\subseteq W(S\setminus R)$ such that
\begin{equation}\label{eq:limit2}
\lim\limits_{n\rightarrow\infty}\left(\rho_n^2\frac{\tilde{G}(\phi_n)+\bar{r}}{2}+\rho_n\tilde{H}_{\bar{r}}(\bar{\rho},\bar{\phi};\phi_n)+\tilde{K}_{\bar{r}}(\bar{\rho},\bar{\phi})\right)=-\infty.
\end{equation}
As in part (a), dropping to a subsequence if necessary we assume $\rho_n\rightarrow\infty$ and $\phi_n\rightarrow\hat{\phi}$ such that $(1,\hat{\phi})\in W(R).$ Note that $(1,\hat{\phi})$ is on the boundary of $W(R).$ Since $(1,\hat{\phi})\in W(R),$ we have $\tilde{G}(\hat{\phi})\geq\tilde{G}.$ In fact, $\tilde{G}(\hat{\phi})>\tilde{G},$ since $\phi\in\tilde{\Phi}\Rightarrow(1,\phi)\in\intt R.$ Hence, $\frac{\tilde{G}(\hat{\phi})+\bar{r}}{2}>0.$ The proof now follows from the same arguments as in Theorem \ref{thm:rsame}.
\item[c)] If $\tilde{G}\leq0$ and $\tilde{\Phi}\neq\tilde{H}_{\bar{r}}^+(\bar{\rho},\bar{\phi})$ then by Theorem \ref{thm:polycone} (c) we have $\bar{x}\not\in\dom e_{\bar{r}}\tilde{f}.$ Since $\dom e_{\bar{r}}f\subseteq\dom e_{\bar{r}}\tilde{f}$ by Lemma \ref{lem:dom1}, we have $\bar{x}\not\in\dom e_{\bar{r}}\tilde{f}.$\qed
\end{itemize}
\begin{rem}
As we saw in Theorem \ref{thm:polycone}, the domain of the Moreau envelope is not identifiable in all situations. For a quadratic function $f$ with general polyhedral domain, we are certain of the domain of $e_{\bar{r}}f$ only in the three situations described in the statement of Theorem \ref{thm:recess}. See Example \ref{ex:badone} for an illustration of how polyhedral domains that are not conic can make it difficult to identify $\dom e_{\bar{r}}f.$
\end{rem}

\subsection{PLQ Functions}

For a quadatic function $f$ whose domain is a single, closed, unbounded polyheral region, we use Theorems \ref{thm:rsame} and \ref{thm:recess} to identify the threshold $\bar{r},$ and $\dom e_{\bar{r}}f.$ We will now use this as a basis for doing the same with a PLQ function. Since a PLQ function is continuous \cite[Proposition 10.21]{rockwets}, every piece is bounded below except possibly those whose domains are unbounded sets. Theorem \ref{thm:main} explicitly identifies the thresholds, and the domains of the Moreau envelopes at the thresholds where possible, of PLQ functions.
\begin{thm}\label{thm:main}
For $i=1,2,\ldots,m,$ let $f_i:\mathbb{R}^n\rightarrow\mathbb{R}$ be quadratic functions on closed, polyhedral domains $S_i:=\dom f_i,$ such that $S_i\cap\intt S_j=\emptyset$ for every $i\neq j,$ and $f_i(x)=f_j(x)$ for all $x\in S_i\cap S_j.$ Let $\bar{r}_i$ be the threshold of $f_i$ for each $i$ (find $\bar{r}_i$ and $\dom e_{\bar{r}_i}f_i$ by applying Theorem \ref{thm:recess} to each $f_i$). Define the function
$$f(x):=\begin{cases}
f_1(x),&x\in S_1,\\
f_2(x),&x\in S_2,\\
&\vdots\\
f_m(x),&x\in S_m.\\
\end{cases}$$
Then the threshold of $f$ is $$\bar{r}=\max\limits_i\{\bar{r}_i\}.$$
Moreover, if we define the active set $\mathcal{A}:=\{i:\bar{r}_i=\bar{r}\},$ then
$$\dom e_{\bar{r}}f=\bigcap\limits_{i\in\mathcal{A}}\dom e_{\bar{r}}f_i.$$
\end{thm}

\textbf{Proof:} We will make use of the following equation in the proof:
\begin{align}
e_rf(\bar{x})&=\inf\limits_{y\in\dom f}\left\{f(y)+\frac{r}{2}\|y-\bar{x}\|^2\right\}\nonumber\\
&=\min\left\{\inf\limits_{y\in S_1}\left\{f_1(y)+\frac{r}{2}\|y-\bar{x}\|^2\right\},\ldots,\inf\limits_{y\in S_m}\left\{f_m(y)+\frac{r}{2}\|y-\bar{x}\|^2\right\}\right\}.\label{eq:mult}
\end{align}
Let $r>\max\limits_i\{\bar{r}_i\}.$ Then by \cite[Theorem 1.25]{rockwets}, we have $e_rf_i(\bar{x})>-\infty$ for all $\bar{x}\in\mathbb{R}^n,$ for all $i.$ Equation \eqref{eq:mult} then gives us that $e_rf(\bar{x})>-\infty$ for all $\bar{x}\in\mathbb{R}^n,$ hence, $\bar{r}\leq\max\limits_i\{\bar{r}_i\}.$ Now let $r<\max\limits_i\{\bar{r}_i\}.$ Then for any $k$ such that $\bar{r}_k=\max\limits_i\{\bar{r}_i\},$ we have $e_rf_k(\bar{x})=-\infty$ for all $\bar{x}\in\mathbb{R}^n.$ Equation \eqref{eq:mult} then gives us that $e_rf(\bar{x})=-\infty$ for all $\bar{x}\in\mathbb{R}^n,$ hence, $\bar{r}\geq\max\limits_i\{\bar{r}_i\}.$ Therefore, $\bar{r}=\max\limits_i\{\bar{r}_i\}.$

If $\bar{r}=0$ and $\dom e_{\bar{r}_i}f_i=\mathbb{R}^n$ for all $i\in\mathcal{A},$ then by Lemma \ref{lem:below} $f_i$ is bounded below for each $i\in\mathcal{A}.$ Since $\max\limits_i\{\bar{r}_i\}=0,$ we know that $\mathcal{A}=\{1,2,\ldots,m\},$ so in fact $f_i$ is bounded below for all $i.$  Hence, $f$ is bounded below as well. By Lemma \ref{lem:below}, $\dom e_{\bar{r}}f=\mathbb{R}^n=\bigcap\limits_{i\in\mathcal{A}}\dom e_{\bar{r}}f_i.$

If we do not have $r=0$ with $\dom e_{\bar{r}_i}f_i=\mathbb{R}^n$ for all $i,$ then consider any $\bar{x}.$ Notice, if $i\notin\mathcal{A},$ then $\bar{r}>\bar{r}_i,$ so $\dom e_{\bar{r}}f_i=\mathbb{R}^n.$ That is, $e_{\bar{r}}f_i(\bar{x})$ is finite. If $i\in\mathcal{A},$ then $e_{\bar{r}}f_i(\bar{x})>-\infty$ if and only if $\bar{x}\in\dom e_{\bar{r}}f_i.$ Hence, we have
$$\dom e_{\bar{r}}f=\bigcap\limits_{i\in\mathcal{A}}\dom e_{\bar{r}}f_i.$$\qed

\begin{rem}Two example applications of Theorem \ref{thm:main} are given in Examples \ref{ex:twoquads} and \ref{ex:six}.\end{rem}

\section{Examples}\label{sec:examples}

We now provide a few examples that illustrate some of the nuances of the results and highlight the procedures given in this paper. The first example illustrates the basic techniques for a full-domain quadratic function.
\begin{ex}\label{ex:basic}
Define $f:\mathbb{R}^2\rightarrow \mathbb{R},$ as the full-domain quadratic
$$f(x):=\frac{1}{2}x^\top\left[\begin{array}{c c}
1&2\\2&-2\end{array}\right]x+\left[\begin{array}{c}1\\1\end{array}\right]x+1.$$
Then the threshold is $\bar{r}=3,$ and $\frac{1}{\bar{r}}\left[\begin{array}{c}1\\1\end{array}\right]\in\dom e_{\bar{r}}f.$
\end{ex}
\textbf{Details:} Let $A=\left[\begin{array}{c c}
1&2\\2&-2\end{array}\right]$  and $b=\left[\begin{array}{c}1\\1\end{array}\right].$ Spectral decomposition of $A$ yields $A=Q^\top DQ$ where $Q=\frac{\sqrt{5}}{5}\left[\begin{array}{c c}2&1\\1&-2\end{array}\right]$ and $D=\left[\begin{array}{c c}2&0\\0&-3\end{array}\right].$ From $D$ we see that $\lambda_1=2$ and $\lambda_2=-3,$ hence $\bar{r}=3.$ As per Theorem \ref{thm:quad3}, we use the substitutions $x=Q^\top y$ and $\bar{x}=Q^\top\bar{y},$ and calculate the Moreau envelope of $f$ at the threshold:
\begin{align*}
e_{\bar{r}}f(Q^\top\bar{y})&=\inf\limits_y\left\{\frac{1}{2}y^\top\left[\begin{array}{c c}2&0\\0&-3\end{array}\right]y+\left(\frac{\sqrt{5}}{5}\left[\begin{array}{c c}2&1\\1&-2\end{array}\right]\left[\begin{array}{c}1\\1\end{array}\right]\right)^\top y+1+\frac{3}{2}\|y-\bar{y}\|^2\right\}\\
&=\inf\limits_y\left\{\left[\frac{5}{2}y_1^2+\left(\frac{3\sqrt{5}}{5}-3\bar{y}_1\right)y_1\right]+\left(-\frac{\sqrt{5}}{5}-3\bar{y}_2\right)y_2+\left(1+\frac{3}{2}(\bar{y}_1^2+\bar{y}_2^2)\right)\right\}\\
&\begin{cases}\frac{3}{5}\bar{y}_1^2-\frac{9\sqrt{5}}{25}\bar{y}_1+\frac{21}{25},&\bar{y}_2=-\frac{\sqrt{5}}{15},\\
-\infty, &\mbox{otherwise.}\end{cases}
\end{align*}
Now we use $\bar{x}=Q^\top\bar{y}$ to find that
$$e_{\bar{r}}f(\bar{x})=\begin{cases}\frac{3}{4}\bar{x}_1^2-\frac{4}{5}\bar{x}_1+\frac{47}{60},&\bar{x}_2=\frac{1}{2}\bar{x}_1+\frac{1}{6},\\-\infty,&\mbox{otherwise.}\end{cases}$$
Hence, we have
$$\dom e_{\bar{r}}f=\left\{\bar{x}:\bar{x}_2=\frac{1}{2}\bar{x}_1+\frac{1}{6}\right\}.$$
Finally, in accordance with Corollary \ref{cor:dom3}, we observe that $\frac{1}{\bar{r}}b\in\dom e_{\bar{r}}f.$\qed

Our next example shows the difficultly in computing $\dom e_{\bar{r}}f$ when non-conic sets are involved.
\begin{ex}\label{ex:badone}
Define $f:\mathbb{R}^2\rightarrow\mathbb{R},$ $f(x,y):=xy.$ Let
\begin{align*}
S_1&=\{(x,y):y=0\}\\
S_2&=\{(x,y):-1\leq y\leq1\},
\end{align*}
and define
\begin{align*}
f_1(x,y)&=\begin{cases}
f(x,y),&(x,y)\in S_1\\
\infty,&\mbox{else,}
\end{cases}\\
f_2(x,y)&=\begin{cases}
f(x,y),&(x,y)\in S_2\\
\infty,&\mbox{else.}
\end{cases}
\end{align*}
Then both $f_1$ and $f_2$ have $G=0$ and $\Phi=H_{\bar{r}}^+(\bar{\rho},\bar{\phi})\setminus H_{\bar{r}}^{++}(\bar{\rho},\bar{\phi}),$ but $\dom e_{\bar{r}}f_1=\mathbb{R}^2,$ whereas $\dom e_{\bar{r}}f_2=\emptyset.$
\end{ex}
\textbf{Details:}
\begin{itemize}
\item[i)] On $S_1,$ the function $f_1$ is equivalently zero. This makes it trivial to find that $G=0$ and $H_{\bar{r}}(\bar{\rho},\bar{\phi};\phi)=0$ for all $x\in\dom f_1,$ for all $\bar{x}\in\mathbb{R}^2.$ Hence, $\Phi=H_{\bar{r}}^+(\bar{\rho},\bar{\phi})\setminus H_{\bar{r}}^{++}(\bar{\rho},\bar{\phi}).$ Since $f_1$ is bounded below, by Lemma \ref{lem:below} we have that $\dom e_{\bar{r}}f_1=\mathbb{R}^2.$
\item[ii)] The recession cone of $S_2$ is $S_1.$ It is left to the reader to verify that $G(\phi)=\sin2\phi,$ $G=\bar{r}=0,$ $\Phi=\{0,\pi\},$ and $H_{\bar{r}}(\bar{\rho},\bar{\phi};\phi)=0,$ so that $\Phi=H_{\bar{r}}^+(\bar{\rho},\bar{\phi})\setminus H_{\bar{r}}^{++}(\bar{\rho},\bar{\phi}).$ Then
\begin{align*}
e_{\bar{r}}f_2(\bar{x},\bar{y})&=\inf\limits_{-1\leq y\leq1}\{xy\}\\
&=-\infty~~~~~~~~~~~~~~~~~\mbox{for all}(\bar{x},\bar{y})\in\mathbb{R}^2.
\end{align*}
Therefore, $\dom e_{\bar{r}}f_2=\emptyset.$\qed
\end{itemize}
Next we have a simple example that shows it possible to construct PLQ functions with equal, positive thresholds, whose Moreau envelope domains are different.
\begin{ex}\label{ex:twoquads}
Define two regions on $\mathbb{R}:$ $S_1=\{x:x\leq0\},~S_2=\{x:x\geq0\}.$ Define
$$\begin{array}{l l c l l}
f_1(x):=-x^2,&x\in S_1,&&f_2(x):=-x^2,&x\in S_2,\\
g_1(x):=-(x+1)^2,&x\in S_1,&&g_2(x):=-(x-1)^2,&x\in S_2.
\end{array}$$
Then the PLQ functions
$$\begin{array}{c c}
F(x):=\begin{cases}
f_1(x),&x\in S_1,\\
f_2(x),&x\in S_2,
\end{cases}&
G(x):=\begin{cases}
g_1(x),&x\in S_1,\\
g_2(x),&x\in S_2,
\end{cases}
\end{array}$$
both have threshold $\bar{r}_f=\bar{r}_g=2,$ but $\dom e_2F=\{0\},$ whereas $\dom e_2G=\emptyset.$
\end{ex}
\begin{figure}[ht]
\begin{center}\includegraphics[scale=0.4]{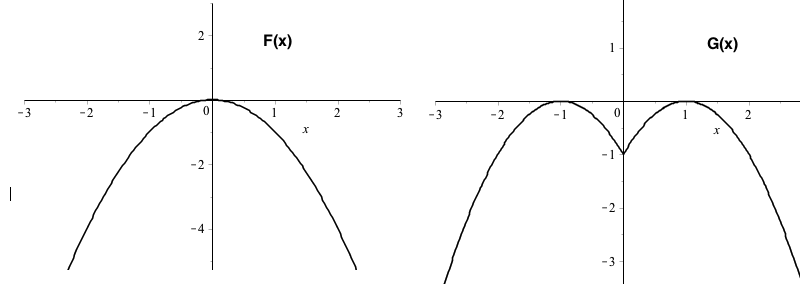}\end{center}
\caption{The Moreau envelopes of PLQ functions with the same $\bar{r}$ may have different domains.}\label{fig:thresh}
\label{fig:rec}
\end{figure}
Figure \ref{fig:rec} makes it easy to see that for $F(x),$ the common real value of the Moreau envelopes is $e_2f_1(0)=e_2f_2(0)=0.$ Hence, $e_2F(0)=0$ and $e_2F(x)=-\infty$ for all $x\neq0,$ which gives $\dom e_2F=\{0\}.$ For $G(x),$ we see that $e_2g_1(-1)=e_2g_2(1)=0$ and the real values of the Moreau envelopes are not at the same point, which gives $e_2G(x)=-\infty$ everywhere. Hence, $\dom e_2G=\emptyset.$\qed

\hbox{}
Finally, we have an example of a six-piece PLQ function on $\mathbb{R}^2.$ We identify the threshold of each piece, and that of the PLQ function. We also make some conclusions with respect to the domain of the Moreau envelope for each piece, and for that of the PLQ function.
\begin{ex}\label{ex:six}
Define six overlapping regions on $\mathbb{R}^2:$
\begin{align*}
S_1=&\{(x,y):y\geq0,~x\leq-2\},\\
S_2=&\{(x,y):x\geq-2,~y\geq x+2,~x\leq0\},\\
S_3=&\{(x,y):y\geq0,~y\leq x+2,~x\leq0\},\\
S_4=&\{(x,y):x\geq0,~y\geq x\},\\
S_5=&\{(x,y):y\geq0,~y\leq x\},\mbox{ and}\\
S_6=&\{(x,y):y\leq0\}.
\end{align*}
Define the quadratic functions
\begin{align*}
f_1(x,y):=&\frac{1}{2}[\begin{array}{c c}x&y\end{array}]\left[\begin{array}{c c}0&-4\\-4&0\end{array}\right]\left[\begin{array}{c}x\\y\end{array}\right]+\left[\begin{array}{c c}1&-3\end{array}\right]\left[\begin{array}{c}x\\y\end{array}\right],\\
f_2(x,y):=&\frac{1}{2}[\begin{array}{c c}x&y\end{array}]\left[\begin{array}{c c}6&-3\\-3&0\end{array}\right]\left[\begin{array}{c}x\\y\end{array}\right]+\left[\begin{array}{c c}7&-1\end{array}\right]\left[\begin{array}{c}x\\y\end{array}\right],\\
f_3(x,y):=&\left[\begin{array}{c c}1&-1\end{array}\right]\left[\begin{array}{c}x\\y\end{array}\right],\\
f_4(x,y):=&\frac{1}{2}[\begin{array}{c c}x&y\end{array}]\left[\begin{array}{c c}12&-7\\-7&0\end{array}\right]\left[\begin{array}{c}x\\y\end{array}\right]+\left[\begin{array}{c c}6&-1\end{array}\right]\left[\begin{array}{c}x\\y\end{array}\right],\\
f_5(x,y):=&\frac{1}{2}[\begin{array}{c c}x&y\end{array}]\left[\begin{array}{c c}0&5\\5&-12\end{array}\right]\left[\begin{array}{c}x\\y\end{array}\right]+\left[\begin{array}{c c}1&4\end{array}\right]\left[\begin{array}{c}x\\y\end{array}\right],\mbox{ and}\\
f_6(x,y):=&\frac{1}{2}[\begin{array}{c c}x&y\end{array}]\left[\begin{array}{c c}0&2\\2&-2\end{array}\right]\left[\begin{array}{c}x\\y\end{array}\right]+[\begin{array}{c c}1&1\end{array}]\left[\begin{array}{c}x\\y\end{array}\right],
\end{align*}
and the PLQ function
$$f(x,y):=\begin{cases}
f_1(x,y),&(x,y)\in S_1,\\
f_2(x,y),&(x,y)\in S_2,\\
f_3(x,y),&(x,y)\in S_3,\\
f_4(x,y),&(x,y)\in S_4,\\
f_5(x,y),&(x,y)\in S_5,\\
f_6(x,y),&(x,y)\in S_6.
\end{cases}$$
Then $f$ has threshold $\bar{r}=\frac{1}{2}+\frac{1}{2}\sqrt{5} \approx 1.618,$ with
    $$\Phi = \{\hat{\phi}\}=\left\{\pi-\arctan \left(\frac{1}{2}+\frac{1}{2}\sqrt{5}\right)\right\}.$$
Moreover, $\dom e_{\bar{r}}f\neq\mathbb{R}^n,~ \dom e_{\bar{r}}f\neq\emptyset.$
\end{ex}

\begin{figure}[ht]
\begin{center}\includegraphics[scale=0.35]{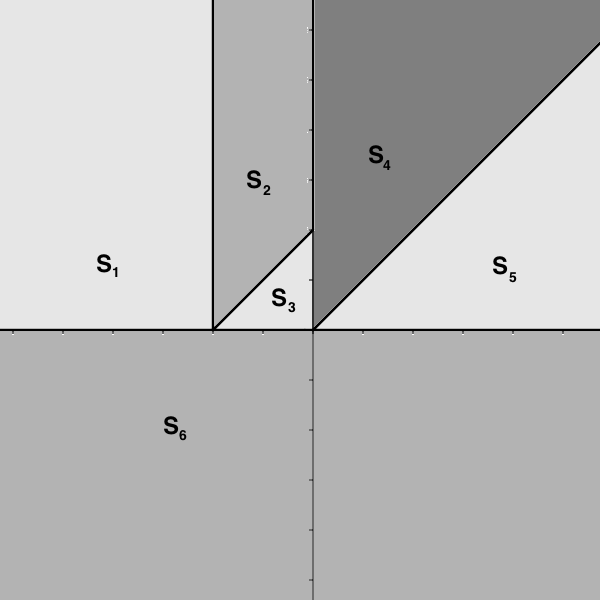}\end{center}
\caption{The partitioning of $\mathbb{R}^2$ for $f(x,y).$}
\label{fig:rec2}
\end{figure}\begin{figure}[ht]
\begin{center}\includegraphics[scale=0.4]{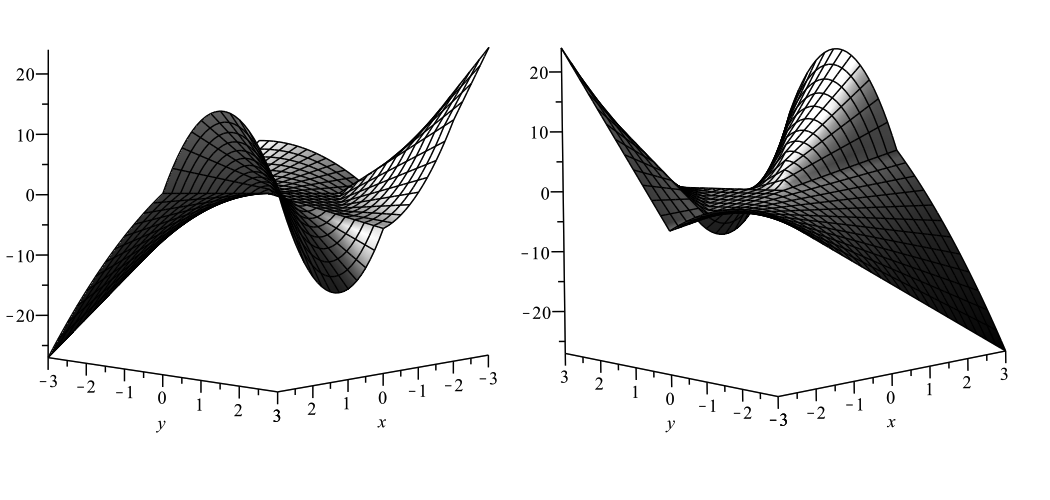}\end{center}
\caption{Two views of the graph of $f(x,y)$.}
\label{fig:rec3}
\end{figure}

\textbf{Details:}
Figure \ref{fig:rec2} shows the six regions of the domain of $f,$ and Figure \ref{fig:rec3} is the graph of $f.$ It is left to the reader to verify that $f$ is indeed a PLQ function, that is, it is continuous at all boundary points.
\begin{itemize}
\item[$S_1:$] This region is not a cone, so we identify the recession cone $R_1$ and use
$$W(R_1)=\left\{(\rho,\phi):\rho\geq0,\phi\in\left[\frac{\pi}{2},\pi\right]\right\}.$$
We consider the restricted function $\tilde{f}_1=f_1$ with $\dom\tilde{f}_1=R_1+(-2,0).$ In polar coordinate form, the function becomes
$$\tilde{f}_1(\rho,\phi)=-4\rho^2\cos\phi\sin\phi+\rho(\cos\phi-3\sin\phi).$$
Then the Moreau envelope at $W((\bar{x},\bar{y}))=(\bar{\rho},\bar{\phi})$ is
$$\inf\limits_{(\rho,\phi)\in W(R_1)}\left\{\rho^2\left(-2\sin2\phi+\frac{r}{2}\right)+\rho\left[\cos\phi-3\sin\phi-r\bar{\rho}\cos(\phi-\bar{\phi})\right]+\frac{r}{2}\bar{\rho}^2\right\}.$$
Using equations (\ref{eq:gphi}) and (\ref{eq:h}), we have $G(\phi)=-2\sin2\phi$ and $H_r(\bar{\rho},\bar{\phi};\phi)=\cos\phi-3\sin\phi-r\bar{\rho}\cos(\phi-\bar{\phi}).$ Notice that $G=\min\limits_{\phi\in\left[\frac{\pi}{2},\pi\right]}G(\phi)=0$ with $\Phi=\argmin\limits_{\phi\in\left[\frac{\pi}{2},\pi\right]}G(\phi)=\left\{\frac{\pi}{2},\pi\right\}.$ This gives $\bar{r}_1=0,$ $H_{\bar{r}_1}\left(\bar{\rho},\bar{\phi};\frac{\pi}{2}\right)=\cos\frac{\pi}{2}-3\sin\frac{\pi}{2}=-3,$ and $H_{\bar{r}_1}(\bar{\rho},\bar{\phi},\pi)=\cos\pi-3\sin\pi=-1,$ independent of the choice of $(\bar{\rho},\bar{\phi}).$ So we have $G\leq0$ and $\Phi\neq H_{\bar{r}_1}^+(\bar{\rho},\bar{\phi})$ for all $\bar{x}\in\mathbb{R}^2.$ Therefore, by Theorem \ref{thm:recess}, $\dom e_{\bar{r}_1}f_1=\emptyset.$
\item[$S_2:$] This region is not a cone, so we identify the recession cone $R_2$ and use
$$W(R_2)=\left\{(\rho,\phi):\rho\geq0,~\phi=\frac{\pi}{2}\right\}.$$
We consider the restricted function $\tilde{f}_2=f_2$ with $\dom\tilde{f}_2=R_2+(-2,0).$ The function in polar coordinates is
$$\tilde{f}_2(\rho,\phi)=3\rho^2\cos^2\phi-3\rho^2\cos\phi\sin\phi+7\rho\cos\phi-\rho\sin\phi,$$
and the Moreau envelope at $(\bar{\rho},\bar{\phi})$ is
$$\inf\limits_{(\rho,\phi)\in W(R_2)}\left\{\rho^2\left(3\cos^2\phi-\frac{3}{2}\sin2\phi+\frac{r}{2}\right)+\rho\left[7\cos\phi-\sin\phi-r\bar{\rho}\cos(\phi-\bar{\phi})\right]+\frac{r}{2}\bar{\rho}^2\right\}.$$
Since we have only one angle in $W(R_2),$ $\phi=\frac{\pi}{2},$ we get $G=3\cos^2\frac{\pi}{2}-\frac{3}{2}\sin\pi=0$ and $\bar{r}_2=0.$ Then $H_{\bar{r}_2}\left(\bar{\rho},\bar{\phi};\frac{\pi}{2}\right)=-1.$ So we have $G\leq0$ and $\Phi\neq H_{\bar{r}_2}^+(\bar{\rho},\bar{\phi})$ for all $\bar{x}\in\mathbb{R}^2.$ Therefore, by Theorem \ref{thm:recess}, $\dom e_{\bar{r}_2}f_2=\emptyset.$
\item[$S_3:$] This region is bounded, so $f_3$ has threshold $\bar{r}_3=0,$ and $\dom e_{\bar{r}_3}f_3=\mathbb{R}^2.$
\item[$S_4:$] This region is a closed, unbounded polyhedral cone, so we use the method of Subsection \ref{sub:gsc}. The function $f_4$ in polar coordinates is
$$f_4(\rho,\phi)=6\rho^2\cos^2\phi-\frac{7}{2}\rho^2\sin2\phi+6\rho\cos\phi-\rho\sin\phi,$$
with domain $W(S_4)=\{(\rho,\phi):\rho\geq0,~\phi\in[\frac{\pi}{4},\frac{\pi}{2}]\}.$ Its Moreau envelope at $(\bar{\rho},\bar{\phi})$ is
$$\inf\limits_{(\rho,\phi)\in W(S_4)}\left\{\rho^2\left(6\cos^2\phi-\frac{7}{2}\sin2\phi+\frac{r}{2}\right)+\rho\left[6\cos\phi-\sin\phi-r\bar{\rho}\cos(\phi-\bar{\phi})\right]+\frac{r}{2}\bar{\rho}^2\right\}.$$
This yields $G(\phi)=6\cos^2\phi-\frac{7}{2}\sin2\phi$ and $G=6\cos^2\hat{\phi}-\frac{7}{2}\sin2\hat{\phi},$ where $\hat{\phi}=\arctan\frac{6+\sqrt{85}}{7}$ is the unique minimizer, hence $G\approx-1.610$ and $\bar{r}_4\approx1.61.$ Since $G<0,$ by Corollary \ref{cor:neq} (noting that $S_4$ is conic) we have $\dom e_{\bar{r}_4}f_4\neq\mathbb{R}^2,$ $\dom e_{\bar{r}_4}f_4\neq\emptyset.$
\item[$S_5:$] This region is also a closed, unbounded polyhedral cone. The function $f_5$ in polar coordinates is
$$f_5(\rho,\phi)=\rho^2(5\cos\phi\sin\phi-6\sin^2\phi)+\rho(\cos\phi+4\sin\phi),$$
with domain $W(S_5)=\{(\rho,\phi):\rho\geq0,~\phi\in[0,\frac{\pi}{4}]\}.$ Its Moreau envelope at $(\bar{\rho},\bar{\phi})$ is
$$\inf\limits_{(\rho,\phi)\in W(S_5)}\left\{\rho^2\left(5\cos\phi\sin\phi-6\sin^2\phi+\frac{r}{2}\right)+\rho\left[\cos\phi+4\sin\phi-r\bar{\rho}\cos(\phi-\bar{\phi})\right]+\frac{r}{2}\bar{\rho}^2\right\}.$$
We find that $G(\phi)$ is minimized uniquely at $\frac{\pi}{4},$ $G=-\frac{1}{2}$ and $\bar{r}_5=\frac{1}{2}.$ Since $G<0,$ by Corollary \ref{cor:neq} we have $\dom e_{\bar{r}_5}f_5\neq\mathbb{R}^2,$ $\dom e_{\bar{r}_5}f_5\neq\emptyset.$
\item[$S_6:$] This region is also a closed, unbounded polyhedral cone. The function $f_6$ in polar coordinates is
$$f_6(\rho,\phi)=\rho^2(2\cos\phi\sin\phi-\sin^2\phi)+\rho(\cos\phi+\sin\phi),$$
with domain $W(S_6)=\{(\rho,\phi):\rho\geq0,~\phi\in[\pi,2\pi]\}.$ Its Moreau envelope at $(\bar{\rho},\bar{\phi})$ is
$$\inf\limits_{(\rho,\phi)\in W(S_6)}\left\{\rho^2\left(2\cos\phi\sin\phi-\sin^2\phi+\frac{r}{2}\right)+\rho\left[\cos\phi+\sin\phi-r\bar{\rho}\cos(\phi-\bar{\phi})\right]+\frac{r}{2}\bar{\rho}^2\right\}.$$
We find that $G(\phi)$ is minimized uniquely at
$$\hat{\phi}=
    \pi-\arctan\frac{10\left[-\frac{1}{200}(50-10\sqrt{5})^{\frac{3}{2}}+\frac{3}{10}\sqrt{50-10\sqrt{5}}\right]}
    {\sqrt{50-10\sqrt{5}}} = \pi - \arctan\left(\frac{1}{2}+\frac{1}{2}\sqrt{5}\right).$$
This provides $G = -\frac{1}{2}-\frac{1}{2}\sqrt{5}$ and $\bar{r}_6 = \frac{1}{2}+\frac{1}{2}\sqrt{5} \approx1.618.$ Since $G<0,$ by Corollary \ref{cor:neq} we have $\dom e_{\bar{r}_6}f_6\neq\mathbb{R}^2,$ $\dom e_{\bar{r}_6}f_6\neq\emptyset.$
\end{itemize}
We summarize these results below. For
$$\hat{\phi}:=\arctan\frac{6+\sqrt{85}}{7},$$
we have the following table.

\begin{table}[ht]\begin{center}
\begin{tabular}{| c | c | c | c|}
\hline
$i$&$r_i$&$r_i$ rounded to $10^{-3}$&$\dom e_{\bar{r}_i}f_i$\\
\hline
1&0&0.000&$\emptyset$\\
2&0&0.000&$\emptyset$\\
3&0&0.000&$\mathbb{R}^2$\\
4&$6\cos^2\hat{\phi}-\frac{7}{2}\sin2\hat{\phi}$&1.610&$\neq\mathbb{R}^2,~\neq\emptyset$\\
5&$\frac{1}{2}$&0.500&$\neq\mathbb{R}^2,~\neq\emptyset$\\
6&$\frac{1}{2}+\frac{1}{2}\sqrt{5}$&1.618&$\neq\mathbb{R}^2,~\neq\emptyset$\\
\hline
\end{tabular}\caption{Results of Example \ref{ex:six}}\label{tab:1}\end{center}\end{table}

By Table \ref{tab:1} and Theorem \ref{thm:main}, $\bar{r}=\bar{r}_6$ and $\dom e_{\bar{r}}f=\dom e_{\bar{r}}f_6.$\qed

\section{Conclusion}\label{sec:conclusion}

In this paper, a variety of methods for identifying the thresholds and domains of Moreau envelopes for functions built on quadratics was presented. Several examples were given to illustrate the techniques. The results found in this paper are applicable to areas of ongoing computational research, wherever calculation of prox-thresholds is needed.

This research raises several questions for further study. For example:
\begin{itemize}
\item[i)] Is it possible to determine computationally the exact threshold of prox-boundedness for some other useful class of functions?
\item[ii)] Any threshold found in this paper, when the domain of the Moreau envelope was the whole space, was equal to zero; does there exist a function $f$ with $\dom e_{\bar{r}}f=\mathbb{R}^n$ such that $\bar{r}>0?$
\item[iii)] Can a calculus of proximal thresholds be created? I.e., given the proximal thresholds of two lsc functions $f$ and $g$, could the proximal thresholds (or bounds for the proximal thresholds) be determined for their sum, product, and composition?
\item[iv)] We relied on the partitioning of $\mathbb{R}^n$ being polyhedral (each region convex, in particular) in order to employ the recession cone for each piece; can this restriction be relaxed?
\item[v)] We also required $n$-dimensional functions, so as to take advantage of the compactness of closed, bounded sets. Can any or all of these results be extended to infinite-dimensional spaces?
\end{itemize}


\begin{thebibliography}{10}

\bibitem{howto}
H.~Bauschke, Y.~Lucet, and M.~Trienis.
\newblock How to transform one convex function continuously into another.
\newblock {\em SIAM Rev.}, 50(1):115--132, 2008.

\bibitem{smoothing}
A.~Beck and M.~Teboulle.
\newblock Smoothing and first order methods: a unified framework.
\newblock {\em SIAM J. Optim.}, 22(2):557--580, 2012.

\bibitem{lawa}
O.~Bretscher.
\newblock {\em Linear Algebra with Application}.
\newblock Prentice-Hall, Upper Saddle River, NJ, 1995.

\bibitem{proxlike}
J.~Chen and S.~Pan.
\newblock A proximal-like algorithm for a class of nonconvex programming.
\newblock {\em Pac. J. Optim.}, 4(2):319--333, 2008.

\bibitem{proximitysums}
P.~Combettes, D.~D{\~u}ng, and B.~V{\~u}.
\newblock Proximity for sums of composite functions.
\newblock {\em J. Math. Anal. Appl.}, 380(2):680--688, 2011.

\bibitem{thresholding}
P.~Combettes and J.~Pesquet.
\newblock Proximal thresholding algorithm for minimization over orthonormal
  bases.
\newblock {\em SIAM J. Optim.}, 18(4):1351--1376, 2007.

\bibitem{linesearch}
R.~Dembo and R.~Anderson.
\newblock An efficient linesearch for convex piecewise-linear/quadratic
  functions.
\newblock In {\em Advances in numerical partial differential equations and
  optimization ({M}\'erida, 1989)}, pages 1--8. SIAM, Philadelphia, PA, 1991.

\bibitem{convhullalg}
B.~Gardiner and Y.~Lucet.
\newblock Convex hull algorithms for piecewise linear-quadratic functions in
  computational convex analysis.
\newblock {\em Set-Valued Var. Anal.}, 18(3-4):467--482, 2010.

\bibitem{ontheconv}
O.~G{\"u}ler.
\newblock On the convergence of the proximal point algorithm for convex
  minimization.
\newblock {\em SIAM J. Control Optim.}, 29(2):403--419, 1991.

\bibitem{proxave}
W.~Hare.
\newblock A proximal average for nonconvex functions: a proximal stability
  perspective.
\newblock {\em SIAM J. Optim.}, 20(2):650--666, 2009.

\bibitem{ncproxave}
W.~Hare and C.~Planiden.
\newblock The {NC}-proximal average for multiple functions.
\newblock {\em Optim. Lett.}, 8(3):849--860, 2014.

\bibitem{parapr}
W.~Hare and C.~Planiden.
\newblock Parametrically prox-regular functions.
\newblock {\em Journal of Convex Analysis}, 21(4), 2014.

\bibitem{proxmap}
W.~Hare and R.~Poliquin.
\newblock Prox-regularity and stability of the proximal mapping.
\newblock {\em J. Convex Anal.}, 14(3):589--606, 2007.

\bibitem{compprox}
W.~Hare and C.~Sagastiz{\'a}bal.
\newblock Computing proximal points of nonconvex functions.
\newblock {\em Math. Program.}, 116(1-2, Ser. B):221--258, 2009.

\bibitem{convprox}
J.~Johnstone, V.~Koch, and Y.~Lucet.
\newblock Convexity of the proximal average.
\newblock {\em J. Optim. Theory Appl.}, 148(1):107--124, 2011.

\bibitem{diffprop}
A.~Jourani, L.~Thibault, and D.~Zagrodny.
\newblock Differential properties of the {M}oreau envelope.
\newblock {\em J. Funct. Anal.}, 266(3):1185--1237, 2014.

\bibitem{bregman}
C.~Kan and W.~Song.
\newblock The {M}oreau envelope function and proximal mapping in the sense of
  the {B}regman distance.
\newblock {\em Nonlinear Anal.}, 75(3):1385--1399, 2012.

\bibitem{ppm}
A.~Kaplan and R.~Tichatschke.
\newblock Proximal point methods and nonconvex optimization.
\newblock {\em J. Global Optim.}, 13(4):389--406, 1998.
\newblock Workshop on Global Optimization (Trier, 1997).

\bibitem{fastmoreau}
Y.~Lucet.
\newblock Fast {M}oreau envelope computation. {I}. {N}umerical algorithms.
\newblock {\em Numer. Algorithms}, 43(3):235--249 (2007), 2006.

\bibitem{whatshape}
Y.~Lucet.
\newblock What shape is your conjugate? {A} survey of computational convex
  analysis and its applications [reprint of mr2496900].
\newblock {\em SIAM Rev.}, 52(3):505--542, 2010.

\bibitem{plqmodel}
Y.~Lucet, H.~Bauschke, and M.~Trienis.
\newblock The piecewise linear-quadratic model for computational convex
  analysis.
\newblock {\em Comput. Optim. Appl.}, 43(1):95--118, 2009.

\bibitem{martreg}
B.~Martinet.
\newblock R\'egularisation d'in\'equations variationnelles par approximations
  successives.
\newblock {\em Rev. Fran\c caise Informat. Recherche Op\'erationnelle}, 4(Ser.
  R-3):154--158, 1970.

\bibitem{moreau1963}
J.-J. Moreau.
\newblock Propri\'et\'es des applications ``prox''.
\newblock {\em C. R. Acad. Sci. Paris}, 256:1069--1071, 1963.

\bibitem{proximite}
J.-J. Moreau.
\newblock Proximit\'e et dualit\'e dans un espace hilbertien.
\newblock {\em Bull. Soc. Math. France}, 93:273--299, 1965.

\bibitem{genhess}
R.~Poliquin and R.~Rockafellar.
\newblock Generalized {H}essian properties of regularized nonsmooth functions.
\newblock {\em SIAM J. Optim.}, 6(4):1121--1137, 1996.

\bibitem{proxfunc}
R.~Poliquin and R.~Rockafellar.
\newblock Prox-regular functions in variational analysis.
\newblock {\em Trans. Amer. Math. Soc.}, 348(5):1805--1838, 1996.

\bibitem{plqopt}
A.~Rantzer and M.~Johansson.
\newblock Piecewise linear quadratic optimal control.
\newblock {\em IEEE Trans. Automat. Control}, 45(4):629--637, 2000.

\bibitem{dynam}
P.A. Rey and C.~Sagastiz{\'a}bal.
\newblock Dynamical adjustment of the prox-parameter in bundle methods.
\newblock {\em Optimization}, 51(2):423--447, 2002.

\bibitem{monops}
R.~Rockafellar.
\newblock Monotone operators and the proximal point algorithm.
\newblock {\em SIAM J. Control Optimization}, 14(5):877--898, 1976.

\bibitem{boundplq}
R.~Rockafellar.
\newblock On the essential boundedness of solutions to problems in piecewise
  linear-quadratic optimal control.
\newblock In {\em Analyse math\'ematique et applications}, pages 437--443.
  Gauthier-Villars, Montrouge, 1988.

\bibitem{rockwets}
R.~Rockafellar and R.~Wets.
\newblock {\em Variational analysis}, volume 317 of {\em Grundlehren der
  Mathematischen Wissenschaften [Fundamental Principles of Mathematical
  Sciences]}.
\newblock Springer-Verlag, Berlin, 1998.

\bibitem{ppa}
W.-Y. Sun, R.J.B. Sampaio, and M.A.B. Candido.
\newblock Proximal point algorithm for minimization of {DC} function.
\newblock {\em J. Comput. Math.}, 21(4):451--462, 2003.

\bibitem{minmoreau}
I.~Yamada, M.~Yukawa, and M.~Yamagishi.
\newblock Minimizing the {M}oreau envelope of nonsmooth convex functions over
  the fixed point set of certain quasi-nonexpansive mappings.
\newblock In {\em Fixed-point algorithms for inverse problems in science and
  engineering}, volume~49 of {\em Springer Optim. Appl.}, pages 345--390.
  Springer, New York, 2011.

\bibitem{funcanal}
K.~Yosida.
\newblock {\em Functional analysis}.
\newblock Die Grundlehren der Mathematischen Wissenschaften, Band 123. Academic
  Press, Inc., New York; Springer-Verlag, Berlin, 1965.
\end{thebibliography}
\end{document}